\documentclass{lms}
\usepackage{amsmath,amssymb}

\let\at@
\catcode`\@=11
\def\atdef@#1{\expandafter\def\csname\space @\string#1\endcsname}
\def\AT#1{\csname\space @\string#1\endcsname}
\let\ampersand@\relax
\newdimen\minaw@
\newdimen\minCDaw@
\def\minCDarrowwidth#1{\RIfMIfI@\onlydmatherr@\minCDarrowwidth
 \else\minCDaw@#1\relax\fi\else\onlydmatherr@\minCDarrowwidth\fi}
\newif\ifCD@
\def\CD{\bgroup\vspace@\relax\let\ampersand@&\iffalse}\fi
 \CD@true\vcenter\bgroup\Let@\tabskip\z@skip\baselineskip20\ex@
 &\hfill$\m@th##$\hfill\crcr}
\def\endCD{\crcr\egroup\egroup\egroup}
\newdimen\bigaw@
\def\setbox@ne{\setbox\@ne}
\def\wd@ne{\wd\@ne}
\def\vspace@{\def\vspace##1{\crcr\noalign{\vskip##1\relax}}}

\catcode`\@\active\let@\AT
\newcount\mscount

\input cd.def
\input dd.def

\newtheorem{l-theorem}[subsubsection]{Theorem}
\newtheorem{l-proposition}[subsubsection]{Proposition}
\newtheorem{l-lemma}[subsubsection]{Lemma}
\newtheorem{l-corollary}[subsubsection]{Corollary}
\newtheorem{l-conjecture}[subsubsection]{Conjecture}
\newnumbered{l-remark}[subsubsection]{Remark}
\newnumbered{l-definition}[subsubsection]{Definition}

\def\theorem{\begin{l-theorem}}
\def\endtheorem{\end{l-theorem}}
\def\proposition{\begin{l-proposition}}
\def\endproposition{\end{l-proposition}}
\def\lemma{\begin{l-lemma}}
\def\endlemma{\end{l-lemma}}
\def\corollary{\begin{l-corollary}}
\def\endcorollary{\end{l-corollary}}
\def\conjecture{\begin{l-conjecture}}
\def\endconjecture{\end{l-conjecture}}
\def\remark{\begin{l-remark}}
\def\endremark{\end{l-remark}}
\def\definition{\begin{l-definition}}
\def\enddefinition{\end{l-definition}}

\def\roster{\begin{enumerate}}
\def\endroster{\end{enumerate}}

\let\local\label
\def\loccit#1{\ref{#1}}
\let\ditto\loccit
\def\iref#1{\ref{#1}\loccit}
\def\itemref#1{\loccit}

\def\period{.\spacefactor1000}

{\catcode`\@=11
\gdef\proclaimfont@{\sl}}

\def\dashes{\begin{itemize}}
\def\enddashes{\end{itemize}}
\def\dash{\item}

\let\AMSbold\mathbf
\let\Bbb\mathbb
\let\Cal\mathcal

\def\bA{\AMSbold A}
\def\bD{\AMSbold D}
\def\bE{\AMSbold E}
\def\bJ{\AMSbold J}
\def\bU{\AMSbold U}
\def\bX{\AMSbold X}

\def\FF{\Bbb F}
\def\DD{\Bbb D}
\def\PP{\Bbb P}

\def\GB{\Bbb B} 
\def\GS{\Bbb S} 
\def\GD{\DD} 

\let\CK\CalK

\let\CG\CalG

\let\CL\CalL
\def\CalD{\Cal D}

\def\tSigma{\tilde\Sigma}
\def\tS{\tilde S}
\def\tX{\tilde X}
\def\tC{\tilde C}
\def\tE{\tilde E}
\def\tQ{\tilde Q}
\def\tU{\tilde U}
\def\tp{\tilde p}
\def\Bb#1{\bar\beta^{#1}}

\def\KK{K}
\def\BK{\bar K}

\let\Ga\alpha
\let\Gb\beta

\def\Tor{\operatorname{Tor}}
\def\discr{\operatorname{discr}}
\def\Pic{\operatorname{Pic}}

\def\Cp#1{\PP^{#1}}
\def\<#1>{\langle#1\rangle}
\def\ie,{\emph{i.e.},}
\def\eg,{\emph{e.g.},}
\def\etc.{\emph{etc}.}

\def\operatorname#1{\mathop{\rm#1}\nolimits}

\def\const{\operatorname{const}}

\def\tr{\operatorname{tr}}
\def\rank{\operatorname{rk}}
\def\Hom{\operatorname{Hom}}
\def\Ker{\operatorname{Ker}}
\def\Ext{\operatorname{Ext}}
\def\Coker{\operatorname{Coker}}
\def\inj{\operatorname{in}}
\def\Tor{\operatorname{Tor}}
\def\Tors{\operatorname{Tors}}
\def\R{\Bbb R}

\def\C{\Bbb C}

\def\Z{\Bbb Z}
\let\ZZ\Z

\let\NN\N
\def\Q{\Bbb Q}
\let\QQ\Q
\let\sminus\smallsetminus

\def\text#1{\hbox{\rm #1}}

\def\qed{\hfill\proofbox}
\let\ge\geqslant
\let\le\leqslant

\def\rom#1{\/{\rm#1}}

\begin{document}

\author{Alex Degtyarev}

\title{Oka's conjecture on irreducible plane sextics}


\classno{%
Primary: 14H30; 
Secondary: 14J28
}

\maketitle

\begin{abstract}
We partially prove and partially disprove Oka's conjecture on the
fundamental group/Alexander polynomial
of an irreducible plane sextic. Among other
results, we enumerate all irreducible sextics with simple
singularities admitting dihedral coverings and find examples
of Alexander equivalent Zariski pairs of irreducible sextics.
\end{abstract}

\section{Introduction}\label{S.intro}

\subsection{Motivation and principal results}
In~\cite{Zariski.group}, O.~Zariski initiated the study
of the fundamental group of the complement of a plane curve as a
topological tool controlling multiple planes ramified at the
curve. He found an example of a curve whose group is not abelian:
it is a sextic with six ordinary cusps which all lie on a conic.
Since then, very few general results have been obtained in this
direction; one may mention M.~V.~Nori's theorem~\cite{Nori},
stating that a curve with sufficiently simple singularities has
abelian fundamental group, and two generalizations of original
Zariski's example, due to B.~G.~Moishezon~\cite{Moishezon} and
M.~Oka~\cite{Oka.Zariski}.

The fundamental group of an algebraic curve~$C$ of large degree is
extremely difficult to compute. As an intermediate tool,
Zariski~\cite{Zariski}
suggested to study its Alexander polynomial $\Delta_C(t)$,
which proved quite
useful in knot theory. This approach was later developed by
A.~Libgober in~\cite{Libgober},~\cite{Libgober2}. The Alexander
polynomial is an algebraic invariant of a group; it is trivial
whenever the group is abelian (see Section~\ref{s.coverings} for
definitions and further references). In the case of plane curves, the
Alexander polynomial can be found in terms of dimensions of
certain linear systems, which depend on the types of the
singular points of the curve and on their global position in~$\Cp2$,
see~\cite{poly}. As a disadvantage, the Alexander polynomial is
often trivial, as it is subject to rather strong divisibility
conditions, see~\cite{Zariski}, \cite{Libgober},
and~\cite{divide}. For example, it is trivial for all irreducible
curves of degree up to five.

The fundamental groups of all curves of degree up to five, both
irreducible and reducible, are known, see~\cite{groups}, and next
degree six has naturally become a subject of intensive research.
A number of contributions has been made
by E.~Artal, J.~Carmona, J.~I.~Cogolludo, C.~Eyral, M.~Oka,
H.~Tokunaga, \etc., see recent survey \cite{Oka.survey}. As a
result, it was discovered that an important r\^ole is played by the so
called sextics of \emph{torus type}, \ie, those whose
equation can be represented in the form $p^3+q^2=0$, where
$p$ and~$q$ are some homogeneous polynomials of degree~$2$
and~$3$, respectively.
Among sextics of torus type is Zariski's
six cuspidal sextic, as well as all other irreducible
sextics with abnormally
large Alexander polynomial, see~\cite{poly}. Furthermore, sextics
of torus type are a principal source of examples of
irreducible curves with nontrivial Alexander polynomial or
nonabelian fundamental group.
Based on the known examples, Oka suggested the
following conjecture.

\conjecture[ (Oka, see~\cite{Oka.conjecture})]\label{conjecture}
Let $C$ be an irreducible plane sextic, which is not of torus type.
Then\rom:
\roster
\item\local{Oka.Delta}
the Alexander polynomial $\Delta_C(t)$ is trivial\rom;
\item\local{Oka.pi.simple}
if all singularities of~$C$ are simple, the group
$\pi_1(\Cp2\sminus C)$ is abelian\rom;
\item\local{Oka.pi}
the fundamental group $\pi_1(\Cp2\sminus C)$ is abelian.
\endroster
\endconjecture


In this paper, we disprove parts~\loccit{Oka.pi.simple}
and~\loccit{Oka.pi} of the conjecture and prove part~\loccit{Oka.Delta}
restricted to sextics with simple singularities (\ie, those of
type~$\bA_p$, $\bD_q$, $\bE_6$, $\bE_7$, or~$\bE_8$,
see~\cite{Arnold} or~\cite{Durfee} for their
definition).

\theorem[ (see Theorem~\ref{torus.equiv} for details)]\label{main.1}
An irreducible plane sextic~$C$ with simple singularities is of
torus type if and only if $\Delta_C(t)\ne1$.
\qed
\endtheorem

\theorem[ (see Theorems~\ref{special.5} and~\ref{special.7} for details)]\label{main.2}
There are irreducible plane sextics~$C_1$, $C_2$ with simple
singularities whose fundamental groups factor to the dihedral
groups~$\GD_{10}$ and~$\GD_{14}$, respectively. The sextics
are not of torus type.
\qed
\endtheorem

\theorem[ (see Theorem~\ref{m-2.groups} for details)]\label{main.3}
There is an irreducible plane sextic~$C$ with a singular point
adjacent to~$\bX_9$ \rom(a quadruple point\rom)
and fundamental group $\GD_{10}\times(\ZZ/3\ZZ)$.
The sextic is not of torus type.
\qed
\endtheorem

Theorems~\ref{main.1}--\ref{main.3} are mere simplified versions of the
statements cited in the titles. We do not prove them separately.

Essentially, Theorem~\ref{main.1} follows from the Riemann-Roch
theorem for $K3$-surfaces, which is not applicable if the curve
has non-simple singular points. For sextics with a singular point
adjacent to~$\bX_9$
(transversal intersection of four smooth branches),
we prove an analog of Theorem~\ref{main.1}
(see Theorem~\ref{m-2.groups}) by calculating the fundamental
groups directly. This result substantiates
Conjecture~\iref{conjecture}{Oka.Delta} in its full version. The
remaining case of curves with a singular point adjacent
to~$\bJ_{10}$
(simple tangency of three smooth branches)
requires a different approach; I am planning to
treat it in a subsequent paper (see~\cite{degt.Oka2}).

\subsection{Other results}
The bulk of the paper is related to the study of irreducible
sextics with simple singularities whose fundamental groups factor
to a dihedral group~$\GD_{2n}$, $n\ge3$. We call such curves
\emph{special}. Alternatively, special is an irreducible sextic
that serves as the ramification locus of a regular
$\GD_{2n}$-covering of the plane.
(Dihedral multiple planes were also extensively studied
by Tokunaga and Artal--Cogolludo--Tokunaga, see recent
papers~\cite{Tokunaga.Galois} and~\cite{ACT}
for further references; their techniques are somewhat more
algebro-geometric.)
We show that only $\GD_6$,
$\GD_{10}$, and $\GD_{14}$ can appear as monodromy groups of
dihedral coverings ramified at irreducible sextics, see
Theorem~\ref{special}, and essentially enumerate all special
sextics, see Sections~\ref{s.statements} and~\ref{s.special}. (The
list of sets of singularities realized by irreducible sextics with
exactly one $\GD_6$-covering is omitted due to its length, and the
rigid isotopy classification of sextics admitting
$\GD_6$-coverings
is not completed. All sets of singularities realized by sextics
of torus type are found in M.~Oka, D.~T.~Pho~\cite{OkaPho.moduli}.)

As a by-product, we discover six sets of singularities that are
realized by both special and non-special irreducible sextics with
$\Delta_C(t)=1$. They give rise to so called
Alexander equivalent Zariski pairs of irreducible sextics
(see Remark~\ref{Zariski.pairs} for details and further
references). To my knowledge, these examples are new. It is worth
mentioning that, as in the case of abundant \emph{vs\period}
non-abundant curves (Zariski pairs of irreducible sextics that
differ by their Alexander polynomials, see~\cite{poly}), within
each pair the special curve is distinguished by the existence of
certain conics passing in a prescribed way through its singular
points. One may hope that, as in the case of abundant curves,
these conics can be used to obtain explicit equations.

The fundamental groups of special sextics are not known. I
would suggest that, at least for the simplest curve in each set,
they are minimal.

\conjecture\label{my.conjecture}
The fundamental groups of the special sextics with the sets of
singularities $3\bA_6$ and $4\bA_4$ are $\GD_{14}\times(\ZZ/3\ZZ)$
and $\GD_{10}\times(\ZZ/3\ZZ)$, respectively.\footnote{Added in
proof: this conjecture has been proved, see~\cite{degt.Oka2},
\cite{degt-Oka}, and~\cite{EyralOka}}
\endconjecture

Any reduced sextic~$C$ of torus type is the critical locus of the
projection to~$\Cp2$ of an irreducible cubic surface
$V\subset\Cp3$. The monodromy of this (irregular) covering is an
epimorphism from $\pi_1(\Cp2\sminus C)$ to the symmetric group
$\GS_3=\GD_6$. Conversely, any such epimorphism gives rise to a
triple covering of~$\Cp2$ ramified at~$C$. We show that the
\emph{existence} of a torus structure is equivalent to the
\emph{existence} of an epimorphism $\pi_1(\Cp2\sminus C)\to\GS_3$,
see Theorem~\ref{torus.equiv}.
(The relation between $\GS_3$-coverings and torus structures was
independently discovered by Tokunaga~\cite{Tokunaga.new};
originally, this question was treated by
Zariski~\cite{Zariski.group}.)
Remarkably, it is not true that
\emph{every} triple plane obtained in this way is a cubic surface.
In the world of irreducible sextics with simple singularities,
there is one counter-example; it is given by Theorem~\ref{9cusps}.

The relation between torus structures and $\GD_6$-coverings is
exploited to detect sextics of torus type and eventually prove
Theorem~\ref{main.1}. Among other results, we classify irreducible
sextics admitting more than one torus structure. The maximal
number is attained at the famous nine cuspidal sextic: it has
twelve torus structures and thirteen $\GD_6$-coverings.

Our study of dihedral coverings is based on
Proposition~\ref{BK2=CK}, which relates the existence of such
coverings to a certain invariant~$\CK_C$ used in the classification of
sextics. As a first step towards reducible curves, we prove
Theorem~\ref{reducible}, which takes into account the $2$-torsion
of the group. Still, this approach can only detect dihedral
quotients of the fundamental group that are compatible with the
standard homomorphism $\pi_1(\Cp2\sminus C)\to\ZZ/2\ZZ$ sending
each van Kampen generator to~$1$. A somewhat complementary
approach was developed by
Tokunaga, see recent paper \cite{Tokunaga.Galois} for further
references. In particular, he constructed a series of
dihedral coverings of the plane ramified at reducible sextics. In
the examples of~\cite{Tokunaga.Galois}, components of the
ramification locus have distinct ramification indices.

Apart from the common goal, Conjecture~\ref{conjecture}, last
section~\ref{S.non-simple} is not related to the rest of the
paper: it is a straightforward application of the results
of~\cite{quintics} and~\cite{groups} dealing with curves of
degree~$m$ with a singular point of multiplicity $m-2$. In
Theorem~\ref{m-2.classes}, we
enumerate all irreducible sextics with a quadruple point and
nonabelian fundamental group. There are seven rigid isotopy
classes; five of them are of torus type, and the remaining two
have trivial Alexander polynomial.

\subsection{Contents of the paper}
In \S\ref{S.tools}, we introduce basic notation and remind a
few facts needed in the sequel. \S\ref{S.sextics} contains a few
auxiliary results, both old and new, related to sextics, Alexander
polynomials, and torus structures. We introduce the notion of
\emph{weight} of a curve, which is used in subsequent statements.
Proposition~\ref{BK2=CK} and
Theorem~\ref{reducible} are also proved here. In
\S\S\ref{S.simple} and~\ref{S.non-simple}, we state and prove
extended versions of Theorems~\ref{main.1}--\ref{main.3}. The most
involved is the case of curves admitting $\GD_6$-coverings.
Technical results obtained in Section~\ref{s.proofs} can be used
in a further study of reducible sextics of torus type.



\section{Preliminaries}\label{S.tools}

\subsection{Basic notation}
For an abelian group~$G$, we use the notation~$G^*$ for the dual group
$\Hom(G,\ZZ)$.
The minimal number of generators of~$G$ is denoted by $\ell(G)$;
if $p$~is a prime, we abbreviate $\ell_p(G)=\ell(G\otimes\FF_p)$.

We use the notation~$\GB_n$ for the braid group on $n$ strings,
$\GS_n$ for the symmetric group of degree~$n$,
and $\GD_{2n}$ for the dihedral group of order~$2n$, \ie,
the semidirect product
$$
1@>>>(\ZZ/n\ZZ)[t]/(t+1)@>>>\GD_{2n}@>>>\ZZ/2\ZZ@>>>1,
$$
where the nontrivial element of $\ZZ/2\ZZ$ acts on the kernel by
the multiplication by~$t$.
One has $\GS_3=\GD_6$. The \emph{reduced braid group} is the
quotient $\GB_n/\Delta^2$ of~$\GB_n$ by its center.
$\GB_3/\Delta^2$ is the free product $(\ZZ/2\ZZ)*(\ZZ/3\ZZ)$.

The Milnor number of an isolated singular point~$P$ is denoted
by~$\mu(P)$. The \emph{Milnor number} $\mu(C)$ of a reduced plane
curve~$C$ is defined as the total Milnor number of all singular
points of~$C$. Given two plane curves~$C$ and~$D$ and an
intersection point $P\in C\cap D$, we use the notation
$(C\cdot D)_P$ for the local intersection index of~$C$ and~$D$
at~$P$.

When a statement is not followed by a proof, either because it is
obvious or because it is cited from another source, it is marked
with
\qed

\subsection{Lattices}\label{s.lattices}
A \emph{lattice} is a finitely generated free abelian group~$L$
equipped with a symmetric bilinear form $b\colon L\otimes L\to\ZZ$.
Usually, we
abbreviate $b(x,y)=x\cdot y$ and $b(x,x)=x^2$. A lattice~$L$ is
\emph{even} if $x^2=0\bmod2$ for all $x\in L$.
As the transition matrix between two integral bases
has determinant $\pm1$, the
determinant $\det L=\det b\in\ZZ$ is well defined.
A lattice~$L$ is called
\emph{nondegenerate} if $\det L\ne0$; it is called
\emph{unimodular} if $\det L=\pm1$.

The bilinear form on a lattice~$L$ extends to $L\otimes\Q$. If
$L$ is nondegenerate, the dual group $L^*$ can be identified
with the subgroup
$$
\bigl\{x\in L\otimes\Q\bigm|
 \text{$x\cdot y\in\ZZ$ for all $x\in L$}\bigr\}.
$$
Hence, $L$ is a subgroup of~$L^*$ and the quotient $L^*\!/L$
is a finite group; it is called the \emph{discriminant group}
of~$L$ and is denoted by $\discr L$. The discriminant
group inherits from $L\otimes\Q$ a symmetric bilinear form
$\discr L\otimes\discr L\to\Q/\ZZ$,
called the \emph{discriminant form},
and, if $L$ is even, its quadratic
extension $\discr L\to\Q/2\ZZ$.
When
speaking about discriminant groups, their
(anti-)isomorphisms, \etc., we assume that the discriminant
form and its quadratic extension are taken
into account. One has
$\mathopen|\discr L\mathclose|=\mathopen|\det L\mathclose|$; in
particular, $\discr L=0$ if and only if $L$ is unimodular.

The orthogonal sum of lattices/discriminant forms is denoted
by~$\oplus$.
Given a lattice~$L$, we use the notation $nL$, $n\in\NN$, for the
orthogonal sum of $n$~copies of~$L$, and $L(q)$, $q\in\QQ$, for
the lattice obtained from~$L$ by multiplying the bilinear form
by~$q$ (assuming that the result is an integral lattice).

From now on, \emph{all lattices are assumed even}.

A \emph{root} in a lattice~$L$ is a vector of
square~$(-2)$. A \emph{root system} is a negative definite lattice
generated by its roots. Each root system admits a unique decomposition
into orthogonal sum of irreducible root systems, the latter
being either $\bA_p$, $p\ge1$, or~$\bD_q$, $q\ge4$,
or~$\bE_6$, $\bE_7$, $\bE_8$. Their discriminant forms are as
follows:
\begin{gather*}
\discr\bA_p=\<-\tfrac{p}{p+1}>,\quad
\discr\bD_{2k+1}=\<-\tfrac{2k+1}4>,\quad
\discr\bD_{2k}=\bmatrix-\frac{k}2&\frac12\\\frac12&1\endbmatrix,\\
\discr\bE_6=\<\tfrac23>,\quad
\discr\bE_7=\<\tfrac12>,\quad
\discr\bE_8=0.
\end{gather*}
Here, $\<p/q>$ with $(p,q)=1$ and $pq=0\bmod2$ represents the
quadratic form on the cyclic group $\ZZ/q\ZZ$
sending~$1$ to $(p/q)\bmod2\ZZ$,
and the $(2\times2)$-matrix represents a quadratic form on the
group $(\ZZ/2\ZZ)^2$.

Recall that the \emph{primitive hull} of a sublattice $L\subset S$
is the sublattice
$$
\tilde L=\bigl\{x\in S\bigm|
 \text{$nx\in L$ for some $n\in\Z$}\bigr\}.
$$
A \emph{finite index extension} of a nondegenerate
lattice~$L$ is a lattice~$S$
containing~$L$ as a finite index subgroup
(\ie, $S$ is the primitive hull of~$L$ in itself),
so that the bilinear
form on~$L$ is the restriction of that on~$S$. Since $S$ is a
lattice, it is canonically embedded into~$L^*$ and the quotient
$\CK=S/L$ is a subgroup of $\discr L$. This subgroup is called
the \emph{kernel} of the extension $S\supset L$. It is
\emph{isotropic}, \ie, the restriction to~$\CK$ of the
discriminant quadratic form is identically zero. Conversely, given
an isotropic subgroup $\CK\subset\discr L$, the group
$S=\{u\in L^*\,|\,(u\bmod L)\in\CK\}$ is a finite index extension
of~$L$.

\theorem[ (see Nikulin~\cite{Nikulin})]
Let $L$ be a nondegenerate even lattice.
Then the map $S\mapsto\CK=S/L\subset\discr L$ establishes a one to one
correspondence between the set of isomorphism classes of finite
index extensions $S\supset L$ and the set of isotropic subgroups
of $\discr L$. Under this correspondence, one has
$\discr S=\CK^\perp\!/\CK$.
\qed
\endtheorem

\subsection{Singularities}
Let $f(x,y)$ be a germ at an isolated singular point~$P$,
let~$\tX$ be the minimal resolution of the singular point~$P$ of
the surface $z^2+f(x,y)=0$, and let $E_i$ be the irreducible
components of the exceptional divisor in~$\tX$.
The group $H_2(\tX)$ is spanned by the classes $e_i=[E_i]$, which
are linearly independent and form a negative definite lattice with
respect to the intersection index form. This lattice is called the
\emph{resolution lattice} of~$P$ and is denoted $\Sigma(P)$.
The basis~$\{e_i\}$
is called
a \emph{standard basis} of~$\Sigma(P)$; it is defined up to
reordering. As usual, $e_i^*$ stand for the elements of the dual
basis of $\Sigma(P)^*=H^2(\tX)$.

\begin{figure}[t]
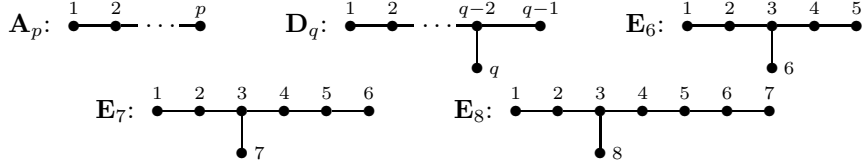

{\normalsize
\begin{gather*}
\DDtop
@[$\bA_p$:\ ]
 @(\bullet)@*<|><1|>@---@(\bullet)@*<|><2|>
 \length7\withdots@---@(\bullet)@*<|><p|>
\endDD
\qquad\quad
\DDtop
@[$\bD_q$:\ ]
 @(\bullet)@*<|><1|>@---@(\bullet)@*<|><2|>
 \length7\withdots@---@(\bullet)@*<|><q-2|>
 \length5@---@(\bullet)@*<|><q-1|>\cr
&&&&&@|---\cr
&&&&&@(\bullet)@*<|q><|>
\endDD
\qquad\quad
\DDtop
@[$\bE_6$:\ ]
 @(\bullet)@*<|><1|>@---
 @(\bullet)@*<|><2|>@---
 @(\bullet)@*<|><3|>@---
 @(\bullet)@*<|><4|>@---
 @(\bullet)@*<|><5|>\cr
&&&&&@|---\cr
&&&&&@(\bullet)@*<|6><|>
\endDD\\\noalign{\smallbreak}
\DDtop
@[$\bE_7$:\ ]
 @(\bullet)@*<|><1|>@---
 @(\bullet)@*<|><2|>@---
 @(\bullet)@*<|><3|>@---
 @(\bullet)@*<|><4|>@---
 @(\bullet)@*<|><5|>@---
 @(\bullet)@*<|><6|>\cr
&&&&&@|---\cr
&&&&&@(\bullet)@*<|7><|>
\endDD
\qquad\quad
\DDtop
@[$\bE_8$:\ ]
 @(\bullet)@*<|><1|>@---
 @(\bullet)@*<|><2|>@---
 @(\bullet)@*<|><3|>@---
 @(\bullet)@*<|><4|>@---
 @(\bullet)@*<|><5|>@---
 @(\bullet)@*<|><6|>@---
 @(\bullet)@*<|><7|>\cr
&&&&&@|---\cr
&&&&&@(\bullet)@*<|8><|>
\endDD
\end{gather*}}
\caption{The standard bases in the Dynkin diagrams}\label{fig.Dynkin}
\end{figure}

If $P$ is simple, of type~$\bA$, $\bD$, $\bE$, then
$\Sigma(P)$ is the irreducible root system of
the same name and one has $\mu(P)=\rank\Sigma(P)$.
In this case, we order the elements of a standard
basis according to the
diagrams shown in Figure~\ref{fig.Dynkin} on
Page~\pageref{fig.Dynkin}.
The order is still defined up to symmetries of the Dynkin graph.

A \emph{rigid isotopy} of plane curves is a topologically
equisingular deformation or, equivalently, a path in a
topologically equisingular stratum of the space of curves. If all
singular points involved are simple, the choice of the category
(topological) in the definition above
is irrelevant, as topologically equivalent simple
singularities are analytically diffeomorphic (see, \eg,~\cite{Arnold}
or~\cite{Durfee}).

\section{Plane sextics}\label{S.sextics}

\subsection{Ramified coverings}\label{s.coverings}
\begin{figure}[t]
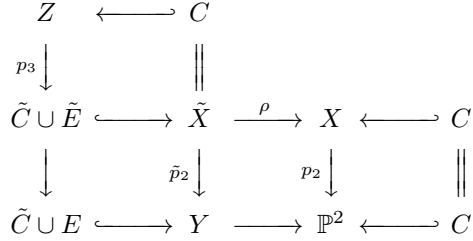

{\normalsize
$$
\CD
@.@.Z@:<---)C\\
@.@.@Vp_3VV@\vert\\
\tC\cup\tE@:(--->\tX@>\rho>>X@:<---)C\\
@VVV@V\tp_2VV@Vp_2VV@\vert\\
\tC\cup E@:(--->Y@>>>\Cp2@:<---)C
\endCD
$$}
\caption{The notation}\label{fig.notation}
\end{figure}

Let $C\subset\Cp2$ be a reduced plane sextic. Throughout the paper
we use the notation introduced in the
diagram in Figure~\ref{fig.notation} on
Page~\pageref{fig.notation}.
Here, $X$ and~$Z$ are, respectively, the double and the $6$-fold
cyclic coverings of~$\Cp2$ ramified at~$C$; clearly, $Z$ can also be
regarded as a triple covering of~$X$. The copies of~$C$ in~$X$
and~$Z$ are identified with~$C$ itself. The map
$\rho\colon \tX\to X$ is the minimal
resolution of singularities of~$X$, $\tC\subset\tX$ is the proper
transform of~$C$, and $\tE$ is the exceptional divisor. The
restriction
$\rho\colon \tX\sminus(\tC\cup\tE)=X\sminus C$
is a diffeomorphism. If all singularities
of~$C$ are simple, then $\tX$ can be obtained as a double covering
of a certain embedded resolution~$Y$ of~$C$; more precisely,
$Y$ is the minimal
resolution in which all odd order components of the pull-back of~$C$ are
smooth and disjoint. The exceptional divisor in~$Y$ is denoted by~$E$.

Let $C_1,\ldots,C_r$ be the irreducible components of~$C$, and let
$\deg C_i=m_i$. From the Poincar\'e duality it follows that the
abelianization of $\pi_1(\Cp2\sminus C)$ is the group
$$\textstyle
H_1(\Cp2\sminus C)=(\ZZ c_1\oplus\ldots\oplus\ZZ c_r)/\sum m_ic_i,
$$
where
$c_i$ is the generator of $H^2(C_i)$ corresponding to the complex
orientation of~$C_i$. The map $c_i\mapsto1$, $i=1,\ldots,r$,
defines canonical epimorphisms
$\pi_1(\Cp2\sminus C)\to\ZZ/6\ZZ\to\ZZ/2\ZZ$. We consider their
kernels
\begin{align*}
\KK_2(C)&=\Ker[\pi_1(\Cp2\sminus C)\to\ZZ/2\ZZ]=
 \pi_1(X\sminus C),\\
\KK_6(C)&=\Ker[\pi_1(\Cp2\sminus C)\to\ZZ/6\ZZ]=
 \pi_1(Z\sminus C)
\end{align*}
and their abelianizations
$$
\BK_2(C)=
 H_1(X\sminus C),\quad
\BK_6(C)=
 H_1(Z\sminus C),
$$
respectively. The deck translations of the coverings~$p_2$ and
$p_2\circ p_3$ induce certain automorphisms $\tr_2$ of $\BK_2(C)$
and $\tr_6$ of $\BK_6(C)$, respectively; the deck translation
of~$p_3$ induces $\tr_6^2$ on $\BK_6(C)$. Group theoretically,
$\tr_2$ and~$\tr_6$ are induced by the conjugation by
(any lifts of)
the generators of $\ZZ/2\ZZ$ and $\ZZ/6\ZZ$, respectively.

The \emph{Alexander polynomial} $\Delta_C(t)$
of a reduced sextic~$C$ can be defined as the characteristic
polynomial of the deck translation automorphism~$\tr_6$ of the
$\C$-vector space $\BK_6(C)\otimes\C=H_1(Z\sminus C;\C)$.
The definition in terms of~$\BK_6$ applies to any group~$G$
equipped with a distinguished epimorphism $G\to\ZZ/6\ZZ$.
One always has
$\Delta_C(t)\mathbin|(t-1)(t^6-1)^4$, see~\cite{Libgober}, and
$\Delta_C(t)$ is defined over~$\QQ$; hence, it is a product of
cyclotomic polynomials.
If $C$ is irreducible, then $\Delta_C(t)\mathbin|(t^2-t+1)^3$,
see~\cite{poly}. Alternative definitions of the Alexander
polynomial of a plane curve and its basic properties can be
found in the original paper~\cite{Libgober} or recent
survey~\cite{Oka.survey}. For the particular case of sextics,
see~\cite{poly} or~\cite{Oka.poly}.

\proposition\label{Delta==>S3}
If $C$ is an irreducible plane sextic with $\Delta_C(t)\ne1$, then
the fundamental group $\pi_1(\Cp2\sminus C)$ factors to the
symmetric group~$\GS_3$.
\endproposition

\proof
Since $\Delta_C(t)\ne1$, the $3$-group $\Hom(\BK_6(C),\FF_3)$
is nontrivial and its order~$3$ automorphism $\tr_6^2$ has
a fixed element.
Hence, $\pi_1(\Cp2\sminus C)$ has a quotient~$G$ which is included
into the exact sequence
$$
1@>>>\ZZ/3\ZZ@>>>G@>>>\ZZ/6\ZZ@>>>1,
$$
so that $\tr_6^2$ acts identically on the kernel. (Here, $\tr_6$
is regarded as a generator of the quotient $\ZZ/6\ZZ$.) Since the
abelianization of $\pi_1(\Cp2\sminus C)$ is $\ZZ/6\ZZ$, the
extension cannot be central. Hence, $\tr_6$ acts on the
kernel
\emph{via}
$x\mapsto-x$,
the exact sequence splits, and $G$ factors to $\GD_6=\GS_3$.
\endproof

\subsection{Sextics of torus type}
A reduced plane sextic~$C$ is said to be of \emph{torus type} if its
equation can be represented in the form
\begin{equation}
p^3(x_0,x_1,x_2)+q^2(x_0,x_1,x_2)=0,
\label{eq.torus}
\end{equation}
where $p$ and~$q$ are some homogeneous polynomials of degree~$2$
and~$3$, respectively. A sextic is of torus type if and only if it
is the critical locus of a projection to~$\Cp2$ of a cubic surface
$V\subset\Cp3$; the latter is given by
$3x_3^3+3x_3p+2q=0$. If $C$ is reduced, then $V$ has isolated
singularities and, hence, is irreducible.

A representation~\eqref{eq.torus}, considered up to scalar
multiples, is called a \emph{torus structure} of~$C$. Each torus
structure gives rise to a conic $Q=\{p=0\}$ and a cubic
$K=\{q=0\}$. The following statement is well known.

\lemma\label{one.torus.structure}
Let $C$ be a plane sextic. Unless $B$ is a union of lines passing
through a single point, each conic $Q$ appears as the conic
$\{p=0\}$ from at most one torus structure of~$C$.
\endlemma

\proof
Let $p^3+q_1^2=\lambda p^3+q_2^2$ be two distinct torus structures
of~$C$. Then $\lambda\ne1$ (as otherwise $q_1=\pm q_2$) and one
has $(\lambda-1)p^3=(q_1-q_2)(q_1+q_2)$. As $\deg p=2$ and
$\deg(q_1\pm q_2)=3$, the conic $Q=\{p=0\}$ must be reducible,
$q=l_1l_2$, and either both $q_1$, $q_2$ belong to the linear
system spanned by $l_1^2l_2$ and $l_1l_2^2$ or
both $q_1$, $q_2$ belong to the linear
system spanned by $l_1^3$ and $l_2^3$. In both cases, $C$ is a
union of lines in the linear system spanned by~$l_1$
and~$l_2$.
\endproof

\remark
Further calculation shows that, in fact, the six lines
constituting an exceptional sextic~$C$ must be in a very special
position. In appropriate coordinates $(x_0:x_1:x_2)$, the two
families are given by
\begin{gather*}
\Ga(x_1x_2)^3+(x_1^3+x_2^3)^2=(\Ga+4)(x_1x_2)^3+(x_1^3-x_2^3)^2,\\
\Ga(x_1x_2)^3+(x_1^2x_2+x_1x_2^2)^2=(\Ga+4)(x_1x_2)^3+(x_1^2x_2-x_1x_2^2)^2.
\end{gather*}
\endremark


Each intersection point $P\in Q\cap K$ is a singular point
for~$C$; such points are called \emph{inner} singularities of~$C$
(with respect to the given torus structures). The other singular
points that~$C$ may have are called \emph{outer}. A simple
calculation using normal forms at~$P$ shows that
an inner singular point can be of type
\dashes
\dash
$\bA_{3k-1}$, if $K$ is nonsingular at~$P$ and $(Q\cdot K)_P=k$,
or
\dash
$\bE_6$, if $K$ is singular at~$P$ and $(Q\cdot K)_P=2$, or
\dash
adjacent to~$\bJ_{10}$ (in the notation of~\cite{Arnold})
otherwise.
\enddashes
(In the latter case, if $P$ is also adjacent to~$\bX_9$ and $B$ is
irreducible, then $P$ is of one of the five last types listed in
Theorem~\ref{m-2.groups}. `Proper' $\bJ$ type singular points are
treated in~\cite{degt.Oka2}.)
Informally, the inner singularities
and their types are due to the topology of the mutual position
of~$Q$ and~$K$, whereas outer singularities occur accidentally in
the family $(\Ga p)^3+(\Gb q)^2=0$ under some special values of
parameters~$\Ga,\Gb\in\C$. A sextic of torus type is called
\emph{tame} if all its singularities are inner. The rigid isotopy
classification of irreducible tame sextics is found
in~\cite{poly}.

\remark
In the case of non-simple points, one should probably speak
about `outer degenerations' of inner singularities. For example,
if $P$ is a node for~$K$ and $(Q\cdot K)_P=3$, the generic inner
singularity at~$P$ is of type $\bJ_{10}=\bJ_{2,0}$. However, under
an appropriate choice of the parameters, it may degenerate
to~$\bJ_{2,1}$ or~$\bJ_{2,2}$. This fact makes the study of sextics
of torus type
with non-simple singularities more involved.
\endremark

\proposition\label{torus==>}
Let $C$ be a reduced sextic of torus type. Then the
group
$\pi_1(\Cp2\sminus C)$ factors to the reduced braid group
$\GB_3/\Delta^2$ and
to the symmetric group~$\GS_3$,
and the Alexander polynomial $\Delta_C(t)$ has at least one factor
$t^2-t+1$.
\endproposition

\proof
All statements follow immediately from the fact that any sextic of
torus type can be perturbed to Zariski's six cuspidal sextic~$C'$,
which is obtained from~$Q$ and~$K$ intersecting transversally at
six points. Hence, there is an epimorphism
$\pi_1(\Cp2\sminus C)\to\pi_1(\Cp2\sminus C')$,
see Zariski~\cite{Zariski.group}.
In the same paper~\cite{Zariski.group}, it is shown that
$\pi_1(\Cp2\sminus C')=\GB_3/\Delta^2$.
\endproof

\definition\label{def.weight}
Let~$P$ be a simple singular point, and let $\Sigma=\Sigma(P)$ be
its resolution lattice. Define the \emph{weight} $w(P)$ as
follows:
$$
w(P)=\begin{cases}
\min\bigl\{-\frac32u^2\bigm|u\in\Sigma^*\sminus\Sigma,\
 3u\in\Sigma\bigr\},\\
0,\quad\text{if $(\discr\Sigma)\otimes\FF_3=0$}.
\end{cases}
$$
The \emph{weight} of a curve~$C$ is the sum of the weights of its
singular points.
\enddefinition

\lemma\label{weight.values}
One has $w(\bA_{3k-1})=k$, $w(\bE_6)=2$, and $w(P)=0$
otherwise. In a standard basis $\{e_i\}$ of~$\Sigma(P)$, the
minimal value of $-\frac32u^2$ as in Definition~\ref{def.weight}
is attained, among other vectors, at~$e_k^*$
or~$e_{2k}^*$ for~$\bA_{3k-1}$ and at~$e_2^*$ or~$e_4^*$
for~$\bE_6$.
\endlemma

\proof
Since $\discr\bE_6=\<\frac23>$, the integer $(3u)^2$ must be
$6\bmod18$. The maximal negative integer with this property
is $-12=(3e_2^*)^2=(3e_4^*)^2$.

Let $\Sigma=\bA_{3k-1}$. Consider the standard representation
of~$\Sigma$ as the orthogonal complement of the characteristic
element $\sum v_i\in\bigoplus\ZZ v_i$, $v_i^2=-1$, $1\le i\le3k$.
An element~$u$ as in
Definition~\ref{def.weight}
has the form
$\frac13\sum m_iv_i$ with
$\sum m_i=0$ and $m_i=\const\ne0\bmod3$ (as all products
$u(v_i-v_j)$ must be integers and $u$ itself must \emph{not} be an
integral vector). If necessary, reversing the sign, we can assume
that $m_i=1\bmod3$.
From the
relation $(m-3)^2+(n+3)^2=(m^2+n^2)-6[(m-n)-3]$ it follows that,
whenever two coefficients in the representation of~$3u$ differ
more than by~$3$, the value of $(3u)^2$ can be increased. Hence,
the coefficients of a square maximizing vector $3u$ take only two values,
which must be $2k$ copies of~$1$ and $k$ copies of~$(-2)$.
Thus, the maximal square is $(3u)^2=-6k$.

For any other irreducible root system~$\Sigma$ one has
$(\discr\Sigma)\otimes\FF_3=0$.
\endproof

\remark\label{weight=d}
From comparing the values given by Lemma~\ref{weight.values}
and those found in~\cite{poly} it follows that,
whenever $w(P)\ne0$, one
has $w(P)=d_{5/6}(P)$, where
$$
d_{5/6}(P)=\#\bigl\{s\in\operatorname{Spec}(P)\bigm|s\le-1/6\bigr\}
$$
are the numbers introduced in~\cite{poly} in
conjunction with
the Alexander polynomial.
(See~\cite{Arnold} for the definition of spectrum.)
Roughly, $d_{5/6}(P)$ is the number of
conditions imposed by~$P$ on the linear system~$\CL_5$ of conics
evaluating~$\Delta_C(t)$.
\endremark

Comparing
Lemma~\ref{weight.values} and
the list of inner
singularities above, one concludes that, for a curve~$C$ of torus
type and conic $Q=\{p=0\}$
defined by a torus structure, $(Q\cdot C)_P=2w(P)$ at each
simple inner
singular point~$P$. (In particular,
if all inner points are simple, then $w(C)\ge\sum_{P\in Q}w(P)=6$.)
The following theorem, which we restate in terms of weights,
asserts that this property is characteristic for conics arising
from torus structures.

\theorem[ (see Degtyarev~\cite{poly} or Tokunaga~\cite{Tokunaga})]\label{torus.criterion}
Let~$C$ be a reduced sextic, and let~$Q$ be a conic \rom(not
necessarily irreducible or reduced\rom)
intersecting~$C$ at simple singular points
so that, at each
intersection
point~$P$, one has $(Q\cdot C)_P=w(P)$.
Then $C$ has a torus structure~\eqref{eq.torus} such that
$Q$ is the conic $\{p=0\}$.
\qed
\endtheorem

The statement proved in~\cite{poly} is stronger
than Theorem~\ref{torus.criterion}: it suffices to require that
the inequality $(Q\cdot C)_P\ge d_{5/6}(P)$ hold
at each intersection point~$P$. In particular, the
intersection points are not restricted \emph{a priori}
to~$\bA_{3k-1}$ or~$\bE_6$.

\subsection{The case of simple singularities}
Let~$C$ be a plane sextic with simple singularities only. Then all
singular points of~$X$ are also simple, and $\tX$ is a
$K3$-surface.
Introduce the following notation:
\dashes
\dash
$L_C=H_2(\tX)$ is the intersection lattice of~$\tX$;
\dash
$h\in L_C$ is the class of the pull-back of
a generic line in~$\Cp2$; one has $h^2=2$;
\dash
$\Sigma_C\subset L_C$ is the sublattice spanned by the
exceptional divisors;
\dash
$S_C=\Sigma_C\oplus\ZZ h$;
\dash
$\tSigma_C$ and~$\tS_C$ are the
primitive hulls of, respectively,~$\Sigma_C$ and~$S_C$ in~$L_C$;
\dash
$\CK_C\subset\discr S_C$ is the kernel of the finite index extension
$\tS_C\supset S_C$.
\enddashes
As is known, $L_C$ is the only unimodular even lattice of signature
$(3,19)$ (one can take $L_C\cong2\bE_8\oplus3\bU$, where $\bU$ is the
hyperbolic plane), and $\Sigma_C$ is the orthogonal sum of the
resolution lattices of all singular points of~$C$.

When a curve~$C$ is understood, we omit subscript~$C$ in the
notation.


The deck translation of the covering $p_2\colon X\to\Cp2$ lifts
to~$\tX$ and permutes the components of~$\tE$; hence, $\tr_2$
induces a certain automorphism of the Dynkin graph of~$\Sigma_C$.
The following lemma is an easy exercise using the embedded
resolution~$Y$ described in Section~\ref{s.coverings}.

\lemma\label{tr.discr}
For each simple singular point~$P$ of~$C$, the automorphism
induced by~$\tr_2$ on the Dynkin graph~$D$ of $\Sigma(P)$ is the
only nontrivial symmetry of~$D$, if $P$ is of type~$\bA_p$,
$\bD_{2k+1}$, or~$\bE_6$
and the identity otherwise. As a consequence, the
induced automorphism of $\discr\Sigma(P)$ is the multiplication
by~$(-1)$.
\qed
\endlemma

The root system~$\Sigma_C$ is called the \emph{set of
singularities} of~$C$. (Since $\Sigma_C$ admits a unique
decomposition into irreducible summands, it does encode the number
and the types of the singular points.)
The triple $h\in S_C\subset L_C$ is called the
\emph{homological type} of~$C$. It is equipped with a natural
orientation~$\theta_C$ of maximal positive definite subspaces in
$S_C^\perp\otimes\R$; it is given by the real and imaginary parts
of the class realized in $H^2(\tX;\C)=L_C\otimes\C$ by a
holomorphic $2$-form on~$\tX$.

\definition\label{def.configuration}
An \emph{\rom(abstract\rom) set of \rom(simple\rom) singularities}
is a root system.
A \emph{configuration} extending a set of
singularities~$\Sigma$ is a finite index
extension $\tS\supset S=\Sigma\oplus\ZZ h$, $h^2=2$,
satisfying the following conditions:
\roster
\item\local{conf.1}
the primitive hull \smash{$\tSigma=h^\perp_{\tS}$}
of~$\Sigma$ in~$\tS$ has no roots other than those in~$\Sigma$;
\item\local{conf.2}
there is no root $r\in\Sigma$ such that $\frac12(r+h)\in\tS$.
\endroster
\enddefinition

\definition\label{def.HT}
An \emph{abstract homological type} extending a set of
singularities~$\Sigma$ is an extension
of $S=\Sigma\oplus\ZZ h$, $h^2=2$, to a lattice
$L\cong2\bE_8\oplus3\bU$ such
that the primitive hull~$\tS$ of~$S$ in~$L$ is a configuration
extending~$\Sigma$.
An abstract homological type is encoded by the triple
$h\in S\subset L$, so that $\ZZ h$ is a direct summand in~$S$ and
$h^\perp_S=\Sigma$.
An \emph{isomorphism} of two abstract
homological types
$h'\in S'\subset L'$ and $h''\in S''\subset L''$ is
an isometry $L'\to L''$
taking~$h'$ to~$h''$ and~$S'$ onto~$S''$.
An
\emph{orientation} of an abstract homological type $h\in S\subset L$ is
an orientation~$\theta$ of maximal positive definite subspaces in
$S^\perp_L\otimes\R$.
\enddefinition

\theorem[ (see Degtyarev~\cite{JAG})]\label{classification}
The homological type $h\in S_C\subset L_C$ of a plane sextic $C$
with simple singularities is an abstract homological type\rom; two
sextics are rigidly isotopic if and only if their oriented
homological types are isomorphic. Conversely, any oriented
abstract homological type is isomorphic to the oriented
homological type of a plane sextic with simple singularities.
\qed
\endtheorem

The existence part of Theorem~\ref{classification}
was first proved by J.-G.~Yang~\cite{Yang}.

\remark\label{rem.Nikulin}
The principal steps of the classification of abstract homological
types are outlined in~\cite{JAG}. A configuration $\tS\supset S$
is determined by its kernel~$\CK$, which plays an important r\^ole
in the sequel. The existence of a primitive extension
$\tS\subset L$ reduces to the existence of an even lattice~$N$ of
signature $(2,20-\rank\tS)$ and discriminant $-\discr\tS$; it can
be detected using Theorem~1.10.1 in~\cite{Nikulin}. Finally, for
the uniqueness one needs to know that
\roster
\item\local{Nik.1}
a lattice $N$ as above is unique up to isomorphism,
\item\local{Nik.2}
each automorphism of $\discr N=-\discr\tS$ is induced by
an isometry of either~$N$ or~$\Sigma$, and
\item\local{Nik.3}
the homological type has an orientation reversing automorphism.
\endroster
In most cases considered in this paper, \loccit{Nik.1} and~\loccit{Nik.2} can
be derived from, respectively, Theorems~1.13.2 and~1.14.2
in~\cite{Nikulin}, and \loccit{Nik.3} follows from the existence of a
vector of square~$2$ in~$N$.
(The case when $N$ is a definite
lattice of rank~$2$ is considered in~\cite{JAG}.) Below, when
dealing with these existence and uniqueness problems,
we just state the result and leave details to the reader.
\endremark

\subsection{Irreducible sextics with simple singularities}\label{s.irreducible}
Our next goal is to relate certain properties of the fundamental
group to the kernel~$\CK_C$ of the extension $\tS_C\supset S_C$.
In this section, we deal with the case of irreducible sextics:
it is more transparent and quite sufficient for the purpose
of this paper. Reducible sextics are considered in
Section~\ref{s.reducible}.

\theorem[(see Degtyarev~\cite{JAG})]\label{irreducible}
A plane sextic~$C$ with simple singularities is irreducible if and
only if the group~$\CK_C$ is free of $2$-torsion.
\qed
\endtheorem

\corollary\label{splitting}
For an irreducible sextic~$C$ with simple singularities one has
$\tS_C=\tSigma_C\oplus\ZZ h$ and $\CK_C=\tSigma_C/\Sigma_C$.
\endcorollary

\proof
One has $\discr(\ZZ h)=\<\frac12>$. Hence, the subgroup
$\CK_C\subset\discr\Sigma_C\oplus\discr(\ZZ h)$ belongs entirely
to $\discr\Sigma_C$, and the orthogonal sum decomposition of~$S_C$
descends to the extension.
\endproof

\corollary\label{inequality}
For an irreducible sextic~$C$ with simple singularities only one has
$\ell_2(\discr\Sigma_C)+\mu(C)\le20$.
\endcorollary

\proof
Since $\CK_C$ is free of $2$-torsion, one has
$\ell_2(\discr\tS)=\ell_2(\discr S)=\ell_2(\Sigma)+1$.
On the other hand, $\ell_2(\discr\tS)\le\rank\tS^\perp=21-\mu(C)$.
\endproof

\proposition\label{BK2=CK}
Let $C$ be an irreducible sextic with simple singularities. Then
$\BK_2(C)$ splits into eigensubgroups,
$\BK_2(C)=\Ker(\tr_2-1)\oplus\Ker(\tr_2+1)$,
and there are
isomorphisms $\Ker(\tr_2-1)=\ZZ/3\ZZ$ and
$\Ker(\tr_2+1)=\Ext(\CK_C,\ZZ)$.
\endproposition

\remark
Proposition~\ref{BK2=CK}, as well as Theorem~\ref{reducible}
below, extend to plane curves of any degree
$(4m+2)$, $m\in\ZZ$: one should
just replace $\ZZ/3\ZZ$ with $\ZZ/(2m+1)\ZZ$ everywhere in the
statements.
\endremark

\proof
One has $\BK_2(C)=H_1(\tX\sminus(\tC\cup\tE))=H^3(\tX,\tC\cup\tE)$
(the Poincar\'e duality) and, since $H^3(\tX)=0$ ($\tX$ is simply
connected), the cohomology exact sequence of pair $(\tX,\tC\cup\tE)$
establishes an isomorphism
$$
\BK_2(C)=\Coker\bigl[\inj^*\colon H^2(\tX)\to H^2(\tC\cup\tE)\bigr].
$$
Let $M=H_2(\tC\cup\tE)$. Since $[\tC]=3h\bmod\Sigma$ in~$L$,
one has $M=\Sigma\oplus(\ZZ\cdot3h)$, the
inclusion homomorphism $\inj_*\colon M\to L$ is
injective,
and the universal coefficients formula implies that
$\BK_2(C)=\Coker[L^*\to M^*]$. The primitive hull of~$M$ in~$L$
is~$\tS$. By definition, $L/\tS$ is a free abelian group;
hence, the adjoint homomorphism
$L^*\to\tS^*$ is onto and one can replace~$L^*$ with~$\tS^*$
in the $\Coker$ expression above. Then, applying
$\Hom(\,\cdot\,,\Z)$ to the free resolution
$0\to M\to\tS\to\tS/M\to0$,
one obtains
$$
\BK_2(C)=\Coker[\tS^*\to M^*]=\Ext(\tS/M,\ZZ).
$$
Due to Corollary~\ref{splitting}, $\tS/M=\CK_C\oplus\ZZ/3\ZZ$ and,
in view of Theorem~\ref{irreducible}, the group $\BK_2(C)$ is free of
$2$-torsion. Hence, $\BK_2(C)$
splits into eigensubgroups of its order~$2$ automorphism~$\tr_2$;
obviously, $h$ is $\tr_2$-invariant,
and the action of~$\tr_2$ on
$\CK_C\subset\discr\Sigma$ is given by Lemma~\ref{tr.discr}.
\endproof

\corollary\label{dihedral}
Let $C$ be an irreducible sextic with simple singularities. Then
there is a canonical one to one correspondence between the set of
normal subgroups $N\subset\pi_1(\Cp2\sminus C)$ with
$\pi_1(\Cp2\sminus C)/N\cong\GD_{2n}$, $n\ge3$,
and the set of subgroups of $\Tor(\CK_C,\ZZ/n\ZZ)$ isomorphic to
$\ZZ/n\ZZ$.
\endcorollary

\proof
The dihedral quotients~$\GD_{2n}$ of the fundamental group are
enumerated by the epimorphisms
$\Ker(\tr_2+1)\to\ZZ/n\ZZ$ modulo
multiplicative units of $(\ZZ/n\ZZ)$,
and the epimorphisms $\Ext(\CK_C,\ZZ)\to\ZZ/n\ZZ$ are the
order~$n$ elements of the group
$$
\Hom(\Ext(\CK_C,\ZZ),\ZZ/n\ZZ)=\Tor(\CK_C,\ZZ/n\ZZ).
$$
(We use the natural isomorphism
$\Hom(\Ext(G,\ZZ),F)=\Tor(G,F)$, which exists for any finite
abelian group~$G$ and any abelian group~$F$.)
\endproof

\subsection{Reducible sextics with simple singularities}\label{s.reducible}
For completeness, we prove an analog of Proposition~\ref{BK2=CK} (and
a more precise version of Theorem~\ref{irreducible}) for reducible
sextics. The results of this section are not used elsewhere in the
paper.

\theorem\label{reducible}
Let $C$ be a reduced plane sextic with simple singularities,
let $C_1,\ldots,C_r$, $r\ge2$, be the irreducible components
of~$C$, and let $c_i\in L=L^*$ be the class realized by the proper
transform~$\tC_i$ of~$C_i$ in~$\tX$, $1\le i\le r$. Then the
following statements hold\rom:
\roster
\item\local{2-torsion}
each residue $c_i\bmod S$ belongs to the subgroup
$\CK_C'=\{\Ga\in\CK_C\,|\,2\Ga=0\}$\rom;
\item\local{2-rank}
the group $\CK_C'$ is generated by the residues $c_i\bmod S$, which
are subject to the only relation $\sum_{i=1}^r c_i=0\bmod S$\rom;
in particular, $\ell_2(\CK_C)=r-1$\rom;
\item\local{torsion}
there is an isomorphism
$\Tors\BK_2(C)=(\ZZ/3\ZZ)\oplus\Ext(\CK_C/\CK_C',\ZZ)$, so that
$\tr_2$ acts \emph{via} $+1$ and~$-1$ on the first and second
summand, respectively\rom;
\item\local{factor}
the group
$\BK_2(C)$ factors to
$(\ZZ/3\ZZ)\oplus\Ext(\CK_C,\ZZ)$, so that
$\tr_2$ acts \emph{via} $+1$ and~$-1$ on the first and second
summand, respectively\rom;
\item\local{free}
the free part $\BK_2(C)/\Tors\BK_2(C)$ is
a free abelian group of rank $r-1$
with the trivial action of~$\tr_2$.
\endroster
\endtheorem

\proof
As in the proof of Proposition~\ref{BK2=CK}, one has a canonical
isomorphism $\BK_2(C)=\Coker[\tS^*\to M^*]$, where
$M=H_2(\tC\cup\tE)$. Now, $M$ is a degenerate lattice, its kernel
being $\Ker[\inj_*\colon M\to L]\cong\ZZ^{r-1}$.
(Indeed, modulo~$\Sigma$ each class~$c_i$ is
homologous to a multiple of~$h$.) This proves
statement~\loccit{free} and gives a natural isomorphism
$\Tors\BK_2(C)=\Ext(\tS/\inj_*M,\ZZ)$, which
reduces~\loccit{torsion} to~\loccit{2-torsion}
and~\loccit{2-rank}.

To prove statement~\loccit{factor}, consider the subgroup
$M_0\subset M$ spanned by the classes of the
exceptional divisors and the total fundamental class
$[\tC]=c_1+\ldots+c_r$. Since the quotient
$M/M_0$ is torsion free, in the diagram
$$
\CD
0@>>>0@>>>\tS^*@=\tS^*@>>>0\\
@.@VVV@VVV@VVV\\
0@>>>(M/M_0)^*@>>>M^*@>>>M_0^*@>>>0
\endCD
$$
the rows are exact, and the $\Ker$--$\Coker$ exact
sequence results in
an epimorphism $\BK_2(C)\to\Ext(\tS/M_0,\ZZ)$. The isomorphism
$\tS/M_0=\ZZ/3\ZZ\oplus\CK_C$ is established similar to
Proposition \ref{BK2=CK}, the two summands being $S/M_0$ and
$\Ker(\tr_2+1)$.

Let $\BK_2(C)\to G$ be the quotient given by~\loccit{factor}. The
further quotient $G/2G$ is an $\FF_2$-vector space
on which $\tr_2$ acts identically.
Hence, $\pi_1(\Cp2\sminus C)$ factors
to an abelian $2$-group~$G'$ with
$\ell_2(G')\ge\dim(G/2G)=\ell_2(\CK_C)$. On the other hand, the
abelianization of~$\pi_1(\Cp2\sminus C)$ is $\ZZ^{r-1}$. Thus,
$\ell_2(\CK_C)\le r-1$.

Let $P$ be a simple singular point, and let
$\Gamma_1,\ldots,\Gamma_s$ be the local branches at~$P$. The
proper pull-back of~$\Gamma_i$ in~$\tX$ represents a certain class
$\gamma_i\in\Sigma(P)^*$, $1\le i\le s$. These classes can easily
be found using the embedded
resolution~$Y$ described in Section~\ref{s.coverings}; it is done
in~\cite{Yang}.
Below, the result is
represented in terms of the basis $\{e_i^*\}$
dual to a standard basis $\{e_i\}$ of $\Sigma(P)$.
(The
representation in terms of the dual basis is very transparent
geometrically: one should just list the exceptional divisors that
intersect the proper transform of a branch.)
\begin{alignat*}4
&\bA_{2k-1}:&\quad&\gamma_1=\gamma_2=e_k^*,
 &\qquad&\bE_6:&\quad&\gamma_1=e_3^*,\\\allowbreak
&\bA_{2k}:&&\gamma_1=e_k^*+e_{k+1}^*,
 &&\bE_7:&&\gamma_1=e_6^*,\ \gamma_2=e_7^*,\\\allowbreak
&\bD_{2k+1}:&&\gamma_1=e_1^*,\ \gamma_2=e_{2k-1}^*,
 &&\bE_8:&&\gamma_1=e_8^*.\\\allowbreak
&\bD_{2k}:&&\gamma_1=e_1^*,\ \gamma_2=e_{2k-1}^*,\ \gamma_3=e_{2k}^*,
\end{alignat*}
On a case by case basis one can verify that
$2\gamma_i=0\bmod\Sigma(P)$, $i=1,\ldots,s$, and
the residues $\gamma_i\bmod\Sigma(P)$ generate the subgroup
$\{\Ga\in\discr\Sigma(P)\,|\,2\Ga=0\}$ and are subject to the only
relation $\sum_{i=1}^s\gamma_i=0\bmod\Sigma(P)$ .

Now, it is obvious that each class~$c_i$, $i=1,\ldots,r$, has the
form
$$
c_i=\frac12(\deg C_i)h+
 \sum\nolimits_{\Gamma_j\subset C_i}\gamma_j,
$$
the sum running over all singular points of~$C$ and all local
branches belonging to~$C_i$. Hence, $2c_i=0\bmod S$. This
proves statement~\loccit{2-torsion} and shows that any nontrivial
relation between the residues $c_i\bmod S$ has the form
$\sum_{i\in I}c_i=0\bmod S$ for some subset
$I\subset\{1,\ldots,r\}$. If both~$I$ and its complement~$\bar I$
are not empty, the curves $C'=\bigcup_{i\in I}C_i$ and
$C''=\bigcup_{i\in\bar I}C_i$ intersect in at least one point~$P$,
which is singular for~$C$. Then, not all local branches at~$P$
belong to~$C'$, and from the properties of classes~$\gamma_j$
stated above it follows that the restriction of
$[C']=\sum_{i\in I}c_i$ to $\Sigma(P)^*$ is not $0\bmod\Sigma(P)$.

Since $r$ residues $c_i\bmod S\in\CK_C'$
are subject to a single relation,
they generate an $\FF_2$-vector space of dimension $r-1$. On the
other hand, as is shown above, $\dim\CK_C'=\ell_2(\CK_C)\le r-1$.
This completes the proof of~\loccit{2-rank} and,
hence,~\loccit{torsion}.
\endproof

\section{Curves with simple singularities}\label{S.simple}

\subsection{Curves of torus type: the statements}\label{s.statements}
In this section, we state our principal results concerning sextics
of torus type. Proofs are given in Section~\ref{s.proofs}.

\theorem\label{torus.equiv}
For an irreducible plane sextic~$C$ with simple singularities, the
following statements are equivalent\rom:
\roster
\item\local{torus.torus}
$C$ is of torus type\rom;
\item\local{torus.Delta}
the Alexander polynomial $\Delta_C(t)$ is nontrivial\rom;
\item\local{torus.B3}
the group $\pi_1(\Cp2\sminus C)$ factors to the
reduced braid group~$\GB_3/\Delta^2$\rom;
\item\local{torus.S3}
the group $\pi_1(\Cp2\sminus C)$ factors to the
symmetric group~$\GS_3$.
\endroster
\endtheorem

A \emph{nine cuspidal sextic} is an irreducible sextic with nine
ordinary cusps, \ie, set of singularities $9\bA_2$. These curves
are well known; they were used by Zariski~\cite{Zariski}
to prove the existence
of non-special six cuspidal sextics. From the Pl\"ucker
formulas it
follows that any nine cuspidal sextic is dual to a nonsingular
cubic curve. In particular, all nine cuspidal sextics
are rigidly isotopic.

\theorem\label{torus.1-1}
Let~$C$ be an irreducible plane sextic with simple singularities,
other than a nine cuspidal sextic. Then
there are canonical bijections between the
following sets\rom:
\roster
\item\local{1.torus}
the set of torus structures on~$C$\rom;
\item\local{1.S3}
the set of normal subgroups $N\subset\pi_1(\Cp2\sminus C)$ with
$\pi_1(\Cp2\sminus C)/N\cong\GS_3$\rom;
\item\local{1.PCK}
the projectivization of the $\FF_3$-vector space
$\CK_C\otimes\FF_3$.
\endroster
In the exceptional case of a nine cuspidal sextic,
still there is a bijection
\loccit{1.S3}~$\leftrightarrow$~\loccit{1.PCK} and an injection
\loccit{1.torus}~$\hookrightarrow$~\loccit{1.PCK}\rom; the image
of the latter injection misses one point.
\endtheorem

The exceptional case in Theorem~\ref{torus.1-1} deserves a
separate statement.

\theorem\label{9cusps}
Let~$C$ be a nine cuspidal sextic. Then there exists one, and only
one, quotient $\pi_1(\Cp2\sminus C)\to\GS_3$ such that the
resulting triple plane $p\colon V\to\Cp2$ ramified at~$C$ is not a cubic
surface. All nine cusps of~$C$ are cusps \rom(Whitney pleats\rom)
of~$p$,
and the covering space~$V$ is a nonsingular surface of
Euler characteristic zero.
\endtheorem

Last three theorems give a detailed description of
sextics of torus
type. In particular, Theorem~\ref{torus.w=8,9} lists all sextics
admitting more than one torus structure. Recall that the
\emph{weight} $w(C)$ of a sextic~$C$ is defined as the total
weight of all its singular points, see
Definition~\ref{def.weight}.

\theorem\label{torus.w}
Let~$C$ be an irreducible sextic with simple singularities.
If the weight $w(C)$ is~$7$ \rom(respectively, $8$ or~$9$\rom),
then $\CK_C=(\ZZ/3\ZZ)^{w(C)-6}$ and $C$ has exactly one
\rom(respectively, four or twelve\rom) torus structures.
If $w(C)=6$, then
$\CK_C=\ZZ/3\ZZ$ or $\CK_C=0$ and $C$ has one or none
torus structure, respectively.
If $w(C)<6$, then $\ell_3(\CK_C)=0$ and $C$ is not of torus type.
\endtheorem

\theorem\label{torus.w=8,9}
Let~$C$ be an irreducible sextic with simple singularities.
If $w(C)=9$, then $C$ is a nine cuspidal sextic.
If $w(C)=8$, then $C$ has one of the following sets of
singularities
\begin{gather*}
8\bA_2,\quad 8\bA_2\oplus\bA_1,\quad
\bA_5\oplus6\bA_2,\quad \bA_5\oplus6\bA_2\oplus\bA_1,\\
2\bA_5\oplus4\bA_2,\quad \bE_6\oplus6\bA_2,\quad
\bE_6\oplus\bA_5\oplus4\bA_2,
\end{gather*}
each set being realized by at least one rigid isotopy
class.\footnote{Added in proof: in fact, each of the eight sets of
singularities in Theorem~\ref{torus.w=8,9} is realized by a single
connected deformation family of sextics, see~\cite{degt.8a2}}
\endtheorem

\theorem\label{torus.w=6}
Let~$C$ be an irreducible sextic with simple singularities and of
weight $w(C)=6$, and assume that $C$ has a singular point of
weight zero other than a simple node \rom(type~$\bA_1$\rom). Then
$\CK_C=\ZZ/3\ZZ$ and $C$ has exactly one torus structure.
\endtheorem

\remark
In the remaining case, $w(C)=6$ and all singular points of weight
zero are simple nodes, the same set of singularities may be
realized by both sextics of torus type and those not of torus
type; they differ by their Alexander polynomials, see abundant
\emph{vs\period} non-abundant curves in~\cite{poly} and~\cite{JAG}.
The first example of this kind is due to
Zariski~\cite{Zariski}.
\endremark

\remark
It is quite straightforward to enumerate all sets of singularities
realized by sextics of torus type; however, the resulting list is
too long.
These sets of singularities are found geometrically in Oka,
Pho~\cite{OkaPho.moduli}.
At present, it is unclear whether
each set of singularities is realized by at most one rigid isotopy
class of sextics of torus type. In many cases, general theorems
of~\cite{Nikulin} do not apply and, considering the amount of
calculations involved, we leave this question open.
\endremark

\subsection{Curves of torus type: the proofs}\label{s.proofs}
Given an integer $w>0$, denote by $\CalD_w$ the direct sum of
$w$~copies of~$\<-\frac23>$; we regard~$\CalD_w$ as an
$\FF_3$-vector space.
Let $\Ga_1,\ldots,\Ga_w$ be some
generators of the summands. An isometry of~$\CalD_w$ is called
\emph{admissible} if it
takes each $\Ga_i$, $i=1,\ldots,w$ to $\pm\Ga_j$ for some~$j$
depending on~$i$.
(One has $\CalD_w=\discr(w\bA_2)$, and the admissible isometries
are those induced by the isometries of $w\bA_2$.)
Define the \emph{weight} $w(\delta)$ of an element
$\delta\in\CalD_w$ as the number of the generators
$\Ga_1,\ldots,\Ga_w$ appearing in~$\delta$ with non-zero
coefficients. Clearly, $\delta$ is isotropic if and only if
$w(\delta)$ is divisible by~$3$.

Let~$C$ be a reduced (not necessarily irreducible) sextic with
simple singularities, and let $w=w(C)$ be the weight of~$C$.
Consider the subgroup $\CG=\CG_C\subset\discr\Sigma_C$ generated
by the elements of order~$3$. Recall that, for each singular
point~$P$ of positive weight $w(P)$, the intersection
$\CG\cap\discr\Sigma(P)$ is generated by a single element $\Gb_P$
of square $-2w(P)/3$. Hence, $\CG$ admits an isometric embedding
to~$\CalD_w$: split the set $G=\{\Ga_1,\ldots,\Ga_w\}$ into
disjoint subsets~$D_P$, assigning $w(P)$ generators to each
singular point~$P$ of positive weight, and map~$\Gb_P$ to
$\sum_{\Ga_i\in D_P}\Ga_i$. Using this embedding,
which is defined up to admissible isometry of~$\CalD_w$,
one can speak about the weights of the elements of~$\CG$.

\lemma\label{star}
In the notation above,
an extension~$\tS$ of the lattice $S=\Sigma\oplus\ZZ h$, $h^2=2$,
defined by an isotropic
subgroup $\CK\subset\CG$ satisfies
condition~\iref{def.configuration}{conf.1} in the definition of
configuration if and only if $\CK$ has the following property\rom:
\roster
\item[($*$)]
each nonzero element of~$\CK$ has weight at least~$6$.
\endroster
\endlemma

\proof
Given $\gamma\in\CK$, the maximal square of a vector
$u\in\tS$ such that $u\bmod S=\gamma$ is $-\frac23w(\gamma)$. This
maximum equals $(-2)$ if and only if $w(\gamma)=3$.
\endproof

\lemma\label{dim*}
Let $w=9$ \rom(respectively, $w=8$ or $w=6,7$\rom), and let
$\CK\subset\CalD_w$ be an isotropic subspace satisfying
condition~\ref{star}$(*)$. Then $\dim\CK\le3$ \rom(respectively,
$\dim\CK\le2$ or $\dim\CK\le1$\rom). Furthermore, a subspace
$\CK_w\subset\CalD_w$
of maximal dimension is unique up to admissible isometry
of~$\CalD_w$\rom; it is generated by
\begin{align*}
w=9\colon &\quad\Ga_1+\ldots+\Ga_9,\quad
      \Ga_1+\Ga_2+\Ga_3-\Ga_4-\Ga_5-\Ga_6,\quad\text{and}\\
     &\quad\Ga_1-\Ga_2+\Ga_4-\Ga_5+\Ga_7-\Ga_8;\\
w=8\colon &\quad\Ga_1+\ldots+\Ga_6\quad\text{and}\quad
           -\Ga_3-\Ga_4+\Ga_5+\ldots+\Ga_8;\\
w\le7\colon &\quad\Ga_1+\ldots+\Ga_6.\\
\end{align*}
\endlemma

\proof
All statements can be proved by a case by case analysis. A more
conceptual proof for the case $w=8$ is given in~\cite{DShapiro},
Lemma~5.2.
This result implies the dimension estimate for $w\le7$ (as the
subgroup of dimension~$2$ involves all eight generators
of~$\CalD_8$) and $w=9$. The uniqueness is obvious in the case
$w\le7$; in the case $w=9$ it can be proved geometrically:
two non-equivalent isotropic subspaces of~$\CalD_9$ of
dimension~$3$ and
satisfying~\ref{star}$(*)$ would give rise to two distinct
configurations extending $9\bA_2$ and,
in view of Theorem~\ref{classification}, to two
rigid isotopy classes of nine cuspidal sextics.
\endproof

\corollary\label{count}
Let $6\le w\le9$, and let $\CK_w\subset\CalD_w$ be the maximal isotropic
subspace given by Lemma~\ref{dim*}. If $w\le8$, then each nonzero
element of~$\CK_w$ has weight~$6$. If $w=9$, then $\CK_w$ has two
elements of weight~$9$ and $24$ elements of weight~$6$.
\qed
\endcorollary

\lemma\label{weight}
Let~$C$ be as above and $w=w(C)$. Assume that the subgroup
$\CG_C\subset\CalD_w$ contains the maximal isotropic
subspace~$\CK_w$ given by Lemma~\ref{dim*}. If $w=8$
\rom(respectively, $w=9$\rom), then for each singular point~$P$
of~$C$ one has $w(P)\le2$ \rom(respectively, $w(P)\le1$\rom).
\endlemma

\proof
Let $G=\bigcup D_P$, $\mathopen|D_P\mathclose|=w(P)$, be the
splitting of the set $G=\{\Ga_1,\ldots,\Ga_w\}$ used to define the
embedding $\CG_C\hookrightarrow\CalD_w$. By the definition
of the embedding, each
element $\gamma=\sum r_i\Ga_i\in\CG_C$ has the following property:
the coefficient function
$c_\gamma\colon i\mapsto r_i$ is constant within each
subset~$D_P$. Since $\CK_w\subset\CG_C$, the maximal weight
$\max w(P)=\max\mathopen|D_P\mathclose|$ is bounded by the maximal
size~$n$ of a subset $D\subset G$ such that the restriction to~$D$
of the coefficient function~$c_\gamma$ of \emph{each} element
$\gamma\in\CK_w$ is a constant (depending on~$\gamma$).
From the description of~$\CK_w$ given in
Lemma~\ref{dim*}, it follows that $n=1$ for $w=9$ and $n=2$ for
$w=8$.
\endproof

\lemma\label{w=6.conic}
In the notation above, there is a natural bijection
between the torus structures of~$C$ and pairs of opposite
elements $\pm\gamma\in\CK_C\cap\CG$ of weight~$6$.
\endlemma

\proof
The statement is essentially contained in~\cite{JAG}, where the
case $w(C)=6$ is considered. Each conic~$Q$ as in
Theorem~\ref{torus.criterion} lifts to two disjoint rational
curves $\tQ_1$, $\tQ_2$ in~$\tX$, and a simple calculation using
the resolution~$Y$ described in Section~\ref{s.coverings}
shows that the fundamental
classes $[\tQ_i]\in L$ have the form
$[\tQ_i]=h+\sum_{P\in Q}\Bb{i}_P$, where,
in a standard basis~$\{e_i\}$ of $\Sigma(P)$,
the elements \smash{$\Bb{i}_P\in\Sigma(P)^*$}, $i=1,2$, are
defined as
\begin{gather*}
\Bb1_P=e_k^*,\quad\Bb2_P=e_{2k}^*\quad
 \text{for $P$ of type $\bA_{3k-1}$},\\
\Bb1_P=e_2^*,\quad\Bb2_P=e_4^*\quad
 \text{for $P$ of type $\bE_6$}.
\end{gather*}
One has
\smash{$(\Bb{i}_P)^2=-\frac23w(P)$}, see
Lemma~\ref{weight.values},
and the residues
\smash{$\Bb{i}_P\bmod\Sigma(P)$}, $i=1,2$, are the two
opposite nontrivial order~$3$ elements of $\discr\Sigma(P)$.
Since one has
$2\sum_{P\in Q}w(P)=C\cdot Q=12$, the residues
$([Q_i]-h)\bmod\Sigma$ form a pair of opposite elements of~$\CK_C$ of
weight~$6$.

Conversely, any order~$3$ element $\gamma\in\CK_C$ can be represented
(possibly, after reordering~$\Bb1$'s and~$\Bb2$'s) as the
residue of the class $\bar\gamma=\sum_{P\in J}\Bb1_P$, the sum
running over a subset~$J$ of the set of singular points with
$\sum_{P\in J}w(P)=w(\gamma)$. If $w(\gamma)=6$, one has
$(\bar\gamma+h)^2=-2$ and, since obviously
$\bar\gamma+h\in\Pic\tX$, the Riemann-Roch theorem implies that
$\bar\gamma+h$ is realized by a (possibly reducible) rational
curve~$\tQ$ in~$\tX$. The image of~$\tQ$ in~$\Cp2$ is a conic~$Q$
as in Theorem~\ref{torus.criterion}.

It remains to notice that, in view of
Lemma~\ref{one.torus.structure}, conics as in
Theorem~\ref{torus.criterion} are in a one-to-one correspondence
with the torus structures of~$C$.
\endproof

\remark
In the proof of Lemma~\ref{w=6.conic}, the lifts~$\tQ_1$, $\tQ_2$
are the connected components of the proper pull-back of~$Q$
provided that $Q$ is nonsingular at each singular point of~$C$. If
$Q$ is singular at a point~$P$ of~$C$, then $P$ is of
type~$\bA_{3k-1}$, $k\ge2$, and a proper pull-back
realizes (locally) a class of the form $e_1^*+e_{k-1}^*$. In this
case, one should include into~\smash{$\tQ_1$} several exceptional
divisors, according to the relation
$e_1^*+e_{k-1}^*+e_1+\ldots+e_{k-1}=e_k^*$. We leave details
to the reader.
\endremark

\lemma\label{kappa}
Let~$C$ be an irreducible sextic with simple singularities, and
let $w(C)\ge7$. Then $\CK_C$ has elements of order~$3$.
\endlemma

\proof
According to~\cite{poly}, the Alexander polynomial $\Delta_C(t)$
is $(t^2-t+1)^s$, where $s$ is the
superabundance of the linear system~$\Cal L_5$
of conics satisfying certain
explicitly described conditions at the singular points of~$C$.
In particular, each singular point~$P$ of positive weight
$w(P)$ imposes $d_{5/6}(P)=w(P)$ conditions, see
Remark~\ref{weight=d}.
Hence,
the virtual dimension of~$\Cal L_5$
is less than $-1$ and
$\Delta_C(t)\ne1$. The statement of the lemma
follows from Proposition~\ref{Delta==>S3} and
Corollary~\ref{dihedral}.
\endproof

\proof[of Theorems~\ref{torus.w} and~\ref{torus.w=8,9}]
First, note that the group~$\CK_C$ has no elements of
order~$9$. Indeed, order~$9$ elements are only present in
$\discr\bA_8$, $2\discr\bA_8$, or~$\discr\bA_{17}$.
However, none of
these discriminants contains an order~$9$ element whose square is
$0\bmod\frac13\ZZ$ (so that it could be compensated by the square
of an order~$3$ element coming from other singular points).

Fix an irreducible sextic~$C$ with simple singularities and
introduce the following notation:
\dashes
\dash
$w=w(C)={}$the weight of~$C$;
\dash
$m={}$the total number of the singular points~$P$ of~$C$ with
$w(P)>0$;
\dash
$e={}$the number of singular points of type~$\bE_6$;
\dash
$\mu'={}$the total Milnor number of the singular points of~$C$ of
weight zero;
\dash
$\kappa=\dim\CK_C\otimes\FF_3$.
\enddashes
The total Milnor number of the singularities of~$C$ is
$\mu=3w-m+e+\mu'$; since $m\le w$, the inequality $\mu\le19$
implies that $w\le9$.

One has $\ell_3(\discr\Sigma)=m$. Hence,
$m-2\kappa\le\ell_3(\discr\tS)\le\rank\tS^\perp=21-\mu$, \ie,
$2\kappa\ge3w+e+\mu'-21$. This inequality, combined with
Lemma~\ref{dim*}, yields:
\dashes
\dash
if $w=9$, then $\kappa=3$ and $e=\mu'=0$;
\dash
if $w=8$, then $\kappa=2$ and $e+\mu'\le1$;
\dash
if $w=7$, then $\kappa=1$ (due to Lemma~\ref{kappa}) and
$e+\mu'\le2$;
\dash
if $w=6$, then $\kappa\le1$ and $e+\mu'\le2\kappa+3$;
\dash
if $w<6$, then $\kappa=0$ (see Lemma~\ref{star}
and~\iref{def.configuration}{conf.1}).
\enddashes
In all cases with $w\ge6$ one has $\mu'\le5$. Furthermore,
whenever $w\ge7$, the subgroup $\CK_C\subset\CG_C\subset\CalD_w$
is the maximal subspace~$\CK_w$ given by Lemma~\ref{dim*}.

To complete the proof of Theorem~\ref{torus.w}, it remains to show
that, whenever $w\ge6$ and $p\ne3$ is a prime, the group~$\CK_C$
is free of $p$-torsion. Since $C$ is irreducible, $\CK_C$ is free
of $2$-torsion, see Theorem~\ref{irreducible}. Thus, we can assume
that $p\ge5$. Then,
$p$-torsion elements are only present in the discriminants
$\discr\bA_i$ with $p\mathbin|(i+1)$. We consider the following two
cases:
\roster
\item
$w(\bA_i)=0$. Then
$i\le\mu'\le5$, leaving the only possibility is $p=5$, $i=4$, \ie,
a single point of type~$\bA_4$.
\item
$w(\bA_i)>0$. Then $3p\mathbin|(i+1)$ and, since $i\le19$, one has
$p=5$, $i=14$. Furthermore, $m\le w-4$
and $e+\mu'=\mu-2w+(m-w)\le3$ (recall that $w\ge6$),
\ie, no other discriminant has elements of prime order $>3$.
\endroster
Thus, one has
$p=5$ and the $5$-torsion of $\discr\Sigma$ comes either from a
single point of type~$\bA_4$ or from a single point of
type~$\bA_{14}$. However, neither $\discr\bA_4$ nor $\discr\bA_{14}$
have an isotropic element of order~$5$.

Prove Theorem~\ref{torus.w=8,9}. If $w(C)=9$, the statement
follows immediately from Lemma~\ref{weight} and the fact that
$\mu'=0$.
If $w(C)=8$, the possible sets of singularities are
easily enumerated using the inequality $e+\mu'\le1$ above
and Lemma~\ref{weight}, which only
allows~$\bA_1$, $\bA_2$, or~$\bE_6$ for a singularity of positive
weight. The sets of
singularities $2\bA_5\oplus4\bA_2\oplus\bA_1$ and
$3\bA_5\oplus2\bA_2$ are ruled out by Corollary~\ref{inequality};
the realizability of the seven sets listed in the
theorem follows from Theorem~\ref{classification} (or
Yang~\cite{Yang}) and
Theorem~1.10.1 in~\cite{Nikulin}, see Remark~\ref{rem.Nikulin}.
\endproof

\proof[of Theorem~\ref{torus.w=6}]
Similar to Lemma~\ref{kappa}, we use the results of~\cite{poly}
(see Remark~\ref{weight=d}) to
evaluate the Alexander polynomial $\Delta_C(t)$. For each singular
point~$P$ other than~$\bA_1$, one has $d_{5/6}(P)\ge1$.
Hence, the total number of conditions on
the conics in~$\CL_5$ is $\sum d_{5/6}(P)\ge w(C)+1=7$. Then, the
virtual dimension of~$\CL_5$ is less than~$-1$, one has
$\Delta_C(t)\ne1$, from Proposition~\ref{Delta==>S3} and
Corollary~\ref{dihedral} it follows that $\CK_C$ has $3$-torsion,
and Theorem~\ref{torus.w} applies.
\endproof

\proof[of Theorem~\ref{torus.equiv}]
The implication \itemref{torus.equiv}{torus.B3}$\implies$\ditto{torus.S3} is
obvious, \ditto{torus.Delta}$\implies$\ditto{torus.S3} is given by
Proposition~\ref{Delta==>S3}, and
\ditto{torus.torus}$\implies$\ditto{torus.B3} and
\ditto{torus.torus}$\implies$\ditto{torus.Delta} are given by
Proposition~\ref{torus==>}. Thus, it remains to show that
\ditto{torus.S3} implies~\ditto{torus.torus}.

Let~$C$ satisfy condition~\ditto{torus.S3}. Due to
Corollary~\ref{dihedral}, the group~$\CK_C$ has elements of
order~$3$ and, comparing Theorem~\ref{torus.w} and
Lemma~\ref{dim*}, one concludes that
$\CK_C\subset\CG_C\subset\CalD_{w(C)}$ is the maximal
isotropic subspace given by Lemma~\ref{dim*}. Then,
Corollary~\ref{count} and Lemma~\ref{w=6.conic} imply that $C$ is
of torus type.
\endproof

\proof[of Theorem~\ref{torus.1-1}]
The bijection \itemref{torus.1-1}{1.S3}~$\leftrightarrow$~\ditto{1.PCK}
is given by Corollary \ref{dihedral}: since
$\CK_C$ has no elements of order~$9$, the order~$3$ subgroups in
$\Tor(\CK_C,\FF_3)$ are in a one to one correspondence with those
in $\CK_C\otimes\FF_3$.

The bijection
\itemref{torus.1-1}{1.torus}~$\leftrightarrow$~\ditto{1.PCK}
is that given by Lemma~\ref{w=6.conic}: in view of
Theorems~\ref{torus.w} and~\ref{torus.w=8,9}, the only exception
is the pair of opposite elements of weight~$9$ that exist in the
case of a nine cuspidal sextic.
\endproof

\proof[of Theorem~\ref{9cusps}]
The triple plane described in the statement corresponds to the two
elements $\pm\gamma\in\CK_C$ of weight~$9$. In general,
let $p\colon V\to\Cp2$ be a triple plane, and assume that the
ramification locus~$C$ has ordinary cusps only. Each cusp~$P$
of~$C$ arises either from a cusp (Whitney pleat) of the
projection~$p$ or from a cusp of~$V$. The former are characterized
by the following property:
\roster
\item"$(*)$"
the composition
$\pi_1(U_P\sminus C)\to\pi_1(\Cp2\sminus C)\to\GS_3$
is an
epimorphism, where $U_P\subset\Cp2$ is a Milnor ball about~$P$.
\endroster

Let~$\tX$, $\tC$, $\tE$ be as in Section~\ref{s.coverings}, and
let $\tU\subset\tX$ be the pull-back of~$U_P$.
Then,
similar to Proposition~\ref{BK2=CK},
one can show that the abelianization of the kernel of
the corresponding homomorphism $\pi_1(U_P\sminus C)\to\ZZ/2\ZZ$
has free part~$\Z$ (with the trivial action of the deck
translation of the covering) and torsion part
$H_1(\partial\tU)=\discr H_2(\tU)=\discr\Sigma(P)$
(see, \eg,~\cite{quartics2} or~\cite{Nikulin}; the discriminant form is
the linking coefficient form on $H_1(\partial\tU)$, which is a
torsion group).
Now, using Proposition~\ref{BK2=CK} and the
identification $\Ext(\CK,\Z)=\Hom(\CK,\Q/\Z)$ (resulting from the
exact sequence $0\to\Z\to\Q\to\Q/\Z\to0$ and the fact that $\CK$ is
finite), one can rewrite the
inclusion homomorphism
$\Tors H_1(\tU\sminus(\tC\cup\tE))\to H_1(\tX\sminus(\tC\cup\tE))$
in the following form:
$$
\discr\Sigma(P)\hookrightarrow\discr\Sigma=
 \Hom(\discr\Sigma,\Q/\Z)\to\Hom(\CK,\Q/\Z).
$$
Here, the first arrow is the inclusion of the direct summand
$\discr\Sigma(P)$, the last one is induced by
the inclusion $\CK\hookrightarrow\discr\Sigma$, and the
isomorphism in the middle is given by the discriminant bilinear
form. In other words, an element $x\in\discr\Sigma(P)$ is sent to
the homomorphism $\CK\to\Q/\Z$, $\gamma\mapsto x\cdot\gamma$.

Thus, a
cusp~$P$ has property~$(*)$ above (for the triple plane $V\to\Cp2$
defined by a pair of opposite elements $\pm\gamma\in\CK$)
if and
only if $\Gb\cdot\gamma\ne0$, where $\Gb$ is a generator of
$\discr\Sigma(P)\cong\Z/3\Z$.
If $w(\gamma)=9$,
this condition holds for
all nine cusps.

The Euler characteristic of~$V$ is found from the Riemann-Hurwitz
formula.
\endproof

\subsection{Other curves admitting dihedral coverings}\label{s.special}
An irreducible sextic is called \emph{special} if its fundamental
group factors to a dihedral group~$\GD_{2n}$, $n\ge3$.
Theorem~\ref{torus.equiv} implies that all irreducible sextics of
torus type are special. In this section, we enumerate other
special sextics with simple singularities.

\theorem\label{special.K}
Let $C$ be an irreducible plane sextic
with simple singularities. Then the group $\Ker(\tr_2+1)$ is
either
$(\ZZ/3\ZZ)^m$, $0\le m\le3$,
or $\ZZ/5\ZZ$, or $\ZZ/7\ZZ$.
\endtheorem

\corollary\label{special}
Let $C$ be an irreducible plane sextic with simple singularities.
Then any dihedral quotient of $\pi_1(\Cp2\sminus C)$ is either
$\GD_6\cong\GS_3$ or $\GD_{10}$ or $\GD_{14}$.
\qed
\endcorollary

\theorem\label{special.7}
There are two rigid isotopy classes of special sextics with simple
singularities whose fundamental group factors to~$\GD_{14}$\rom;
their sets of singularities are $3\bA_6$ and $3\bA_6\oplus\bA_1$.
The set of singularities~$3\bA_6$ can also be realized by a
non-special irreducible sextic.

The two special sextics above can be characterized as follows\rom:
there is an ordering $P_1$, $P_2$, $P_3$ of the three $\bA_6$
points such that, for every cyclic permutation $(i_1i_2i_3)$, there
is a conic whose local intersection index with~$C$ at~$P_{i_k}$
equals $2k$.
\endtheorem

\theorem\label{special.5}
There are eight rigid isotopy classes of special sextics with simple
singularities whose fundamental group factors to~$\GD_{10}$\rom;
each class is determined by its set of singularities, which is one
of the following\rom:
\begin{gather*}
4\bA_4,\quad 4\bA_4\oplus\bA_1,\quad 4\bA_4\oplus2\bA_1,\quad
4\bA_4\oplus\bA_2,\\
\bA_9\oplus2\bA_4,\quad \bA_9\oplus2\bA_4\oplus\bA_1,\quad
\bA_9\oplus2\bA_4\oplus\bA_2,\quad
2\bA_9.
\end{gather*}
The sets of singularities
$4\bA_4$, $4\bA_4\oplus\bA_1$,
$\bA_9\oplus2\bA_4$, $\bA_9\oplus2\bA_4\oplus\bA_1$, and
$2\bA_9$ are also realized by non-special irreducible sextics.

The eight special sextics above can be characterized as
follows\rom: there are two conics~$Q_1$, $Q_2$ with the following
properties\rom:
\dashes
\dash
$Q_1$ and~$Q_2$ intersect transversally at each singular point
of~$C$ of type~$\bA_4$, and they have a simple tangency at each
singular point of~$C$ of type~$\bA_9$\rom;
\dash
at each singular point of type~$\bA_4$, the local intersection
indices of~$C$ with the two conics are~$2$ and~$4$\rom;
\dash
at each singular point of type~$\bA_9$, the local intersection
indices of~$C$ with the two conics are~$4$ and~$8$.
\enddashes
\endtheorem

\remark\label{Zariski.pairs}
The Alexander polynomials of all curves listed in
Theorems~\ref{special.7} and~\ref{special.5} are trivial, \eg, due
to Proposition~\ref{Delta==>S3}.
Hence, the sets of singularities $3\bA_6$, $4\bA_4$, $4\bA_4\oplus\bA_1$,
$\bA_9\oplus2\bA_4$, $\bA_9\oplus2\bA_4\oplus\bA_1$, and
$2\bA_9$ that are realized by both special and non-special curves
give rise to Alexander equivalent Zariski pairs of
\emph{irreducible} sextics. This means that two irreducible
curves~$C_1$, $C_2$ share the same set of singularities
and Alexander polynomial but have non-diffeomorphic complements
$\Cp2\sminus C_i$ (see~\cite{Artal.def} for precise definitions).
In our case, the fundamental groups $\pi_1(\Cp2\sminus C)$ differ:
one does and the other does not admit dihedral quotients.
\endremark

\proof[of Theorem~\ref{special.K}]
Since $C$ is irreducible, $\CK_C$ is free of $2$-torsion.
The case when $\CK_C$ has $3$-torsion is considered in
Theorem~\ref{torus.w}. For a prime $p\ge5$, any simple singularity
whose discriminant has elements of order $p^a$ is of type $\bA_i$
with $p^a\mathbin|(i+1)$. Since the total Milnor number $\mu\le19$,
one has $p^a=5$, $7$, $11$, $13$, $17$, or~$19$. In the last four
cases, $p=11$, $13$, $17$, or~$19$, the set of singularities has
at most one
point with $p$-torsion in the discriminant, which is of type
$\bA_{p-1}$; however, $\discr\bA_{p-1}=\<-\frac{p-1}p>$ does not
have isotropic elements of order~$p$. The remaining cases $p=7$ and $p=5$
are considered in Theorems~\ref{special.7} and~\ref{special.5},
respectively. In particular,
it is shown that the $p$-primary part of~$\CK_C$ is $\ZZ/p\ZZ$.
Comparing the sets of singularities listed in
Theorems~\ref{special.7} and~\ref{special.5}, one immediately
concludes that $\CK_C$ cannot have both $7$- and $5$-torsion.
\endproof

\proof[of Theorem~\ref{special.7}]
Since~$\mu\le19$, the part of~$\Sigma$ whose discriminant has
$7$-torsion is either $a\bA_6$, $1\le a\le3$, or
$\bA_{13}\oplus\bA_6$.
(As explained in Section~\ref{s.lattices}, the discriminants of
irreducible root systems of type~$\bD$ and~$\bE$ are $2$-groups.)
It is easy to see that only
$\discr(3\bA_6)$ contains an order~$7$ isotropic element,
and it
is unique up to isometry of~$\Sigma$.
Besides, since $\discr(3\bA_6)=(\ZZ/7\ZZ)^3$ and the form is
nondegenerate, this group cannot contain an isotropic subgroup
larger than $\ZZ/7\ZZ$. These observations restrict the possible
sets of singularities to those listed in the statement.
The existence of all three
curves mentioned in the statement and the uniqueness of the two
special curves are straightforward, see
Theorem~\ref{classification} and Remark~\ref{rem.Nikulin};
the set of singularities $3\bA_6\oplus\bA_1$ cannot be realized
by a non-special curve since for such a curve one would have
$\ell_7(\discr\tS)=3>2=\rank\tS^\perp$.

The characterization of the special curves in terms of conics is
obtained similar to Lemma~\ref{w=6.conic}.
Let $P_i$, $i=1,2,3$, be the three points of type~$\bA_6$, and
denote by $e_{ij}$, $j=1,\ldots,6$, a standard basis of
$\Sigma(P_i)$. Then, up to a symmetry of the Dynkin graph,
an isotropic element in $\discr\Sigma$ is given by
$\gamma=e_{11}^*+e_{22}^*+e_{33}^*\bmod\Sigma$. The class
$e_{11}^*+e_{22}^*+e_{33}^*+h$ has square~$(-2)$; hence, it is
realized by a rational curve in~$\tX$, which projects to a conic
in~$\Cp2$. The two other conics are obtained from
the classes
$2\gamma=e_{12}^*+e_{24}^*+e_{36}^*\bmod\Sigma$ and
$3\gamma=e_{13}^*+e_{26}^*+e_{32}^*\bmod\Sigma$.
Conversely, each conic as in the statement lifts to a rational
curve in~$\tX$ that realizes an order~$7$ element in
$\discr\Sigma$.
\endproof

\proof[of Theorem~\ref{special.5}]
The singularities whose discriminants contain elements of
order~$5$ are $\bA_4$, $\bA_9$, $\bA_{14}$, and~$\bA_{19}$. The
only imprimitive finite index extension of $2\bA_4$ is
$2\bA_4\subset\bE_8$; it violates
condition~\iref{def.configuration}{conf.1} in the definition of
configuration. With this possibility ruled out, the discriminants
containing an order~$5$ isotropic subgroup are those of $4\bA_4$,
$\bA_9\oplus2\bA_4$, and $2\bA_9$. In each case, the subgroup is
unique up to a symmetry of the Dynkin graph; it is generated by
the residue $\gamma=\bar\gamma\bmod\Sigma$, where
$\bar\gamma$ is given by
$$
e_{11}^*+e_{21}^*+e_{32}^*+e_{42}^*,\quad
f_{14}^*+e_{11}^*+e_{21}^*,\quad\text{and}\quad
f_{14}^*+f_{22}^*,
$$
respectively. Here, $\{e_{ij}\}$, $j=1,\ldots,4$, is a standard
basis in the $i$-th copy of~$\bA_4$, and $\{f_{kj}\}$,
$j=1,\ldots,9$, is a standard basis in the $k$-th copy of~$\bA_9$.
Similar to Lemma~\ref{w=6.conic}, these expressions give a
characterization of the special curves in terms of conics: the class
$\bar\gamma+h$ has square~$(-2)$ and is realized by a rational
curve in~$\tX$; its projection to~$\Cp2$ is one of the two conics.
The other conic is obtained from a similar representation
of~$2\gamma$.

The rest of the theorem is
an application of Theorem~\ref{classification} and Nikulin's
results on lattices.
The sets of singularities $4\bA_4\oplus3\bA_1$,
$\bA_9\oplus2\bA_4\oplus2\bA_1$, and
$2\bA_9\oplus\bA_1$ cannot be realized by irreducible curves due
to the genus formula
(alternatively, due to
Corollary~\ref{inequality}).
The sets of
singularities $4\bA_4\oplus\bA_2\oplus\bA_1$ and
$4\bA_4\oplus\bA_3$ do not extend to abstract homological types
due to Theorem~1.10.1 in~\cite{Nikulin}, see
Remark~\ref{rem.Nikulin}.
The (non-)existence of the other curves mentioned in the statement is
given by Theorem~\ref{classification}, see
Remark~\ref{rem.Nikulin}. The uniqueness of the special curves is
also given by Theorem~\ref{classification}: in most cases one can
apply either Theorem~1.14.2 in~\cite{Nikulin} or Theorem~1.13.2
in~\cite{Nikulin} and the fact that all automorphisms of
$\discr S$ are realized by
isometries of~$\Sigma$. (The
last statement is true in all cases except
$\bA_9\oplus2\bA_4\oplus\bA_1$.) We leave details to the
reader.
\endproof

\section{Some curves with a non-simple singular point}\label{S.non-simple}

In this concluding section we try to substantiate
Conjecture~\iref{conjecture}{Oka.Delta} extended to all irreducible
sextics. Here, we consider sextics with a non-simple singular
point adjacent to~$\bX_9$, \ie, a point of multiplicity~$4$
or~$5$. The only remaining case of a singular point adjacent
to~$\bJ_{10}$ will be dealt with in a separate paper.

\subsection{Sextics with a singular point of multiplicity~$5$}\label{s.m5}
This case is trivial: due to~\cite{groups}, any irreducible sextic
with a singular point of multiplicity~$5$ has abelian fundamental
group and, hence, trivial Alexander polynomial. Note that
Proposition~\ref{torus==>} implies that none of such sextics is of
torus type.

\subsection{Sextics with a singular point of multiplicity~$4$}\label{s.m4}
The rigid isotopy classification of plane curves~$C$
with a singular point~$P$ of multiplicity $\deg C-2$ is found
in~\cite{quintics}. Let $m=\deg C$.
In appropriate coordinates $(x_0:x_1:x_2)$
the curve is given by a polynomial of the form
$$
x_0^2a(x_1,x_2)+x_0b(x_1,x_2)+c(x_1,x_2),
$$
where $a$, $b$, and~$c$ are some homogeneous
polynomials of degree $m-2$, $m-1$, and~$m$, respectively. The
discriminant $D=b^2-4ac$ has degree~$2m-2$. (It is required that
$D$ is not identically zero.) Since $C$ is assumed irreducible,
$a$, $b$, and~$c$ should not have common roots. Let $x_i$,
$i=1,\ldots,k$, be all distinct roots of~$aD$.
The \emph{formula}
of~$C$ is defined as the (unordered) set $\{(p_i,q_i)\}$,
$i=1,\ldots,k$, where $p_i$ and~$q_i$ are the multiplicities
of~$x_i$ in~$a$ and~$D$, respectively. The formula of an
irreducible curve of degree~$m$ has the following properties:
\roster
\item\local{m4.1}
$\sum_{i=1}^kp_i=m-2$, and $\sum_{i=1}^kq_i=2m-2$;
\item\local{m4.2}
for each~$i$, either $p_i=q_i$ or the smallest of~$p_i$, $q_i$ is
even;
\item\local{m4.3}
at least one of~$q_i$ is odd.
\endroster
An \emph{elementary equivalence} of a formula is replacing two
pairs $(1,0)$, $(0,1)$ with one pair $(1,1)$. Geometrically, this
procedure means that a `vertical' tangency point of~$C$ disappears
at infinity making~$P$ an inflection point of one of its smooth
branches. Clearly, this is a rigid isotopy.

\theorem[ (see Degtyarev~\cite{quintics})]\label{m-2.classes}
Two irreducible curves of degree~$m$, each with a singular point
of multiplicity~$m-2$, are rigidly isotopic if and only if their
formulas are related by a sequence of elementary equivalences and
their inverses. Any set of pairs of nonnegative integers
satisfying conditions~\itemref{s.m4}{m4.1}--\ditto{m4.3} above is realized
as the formula of an irreducible curve of degree~$m$.
\qed
\endtheorem

Let~$G_p$, $p=2,4$, be the group given by
$$
G_p=\bigl<u,v\bigm|u^p=v^p,\ (uv)^2u=v(uv)^2,\ v^4=(uv)^5\bigr>.
$$
In~\cite{groups}, it is shown that $\mathopen|G_4\mathclose|=\infty$
and $\mathopen|G_2\mathclose|=30$;
there is a split exact sequence
$$
1@>>>\FF_5[t]/(t+1)@>>>G_2@>>>\ZZ/6\ZZ@>>>1,
$$
$t$ being the conjugation action on the kernel of a generator of
$\ZZ/6\ZZ$. The Alexander polynomials of both~$G_2$ and~$G_4$ are
trivial.

\theorem\label{m-2.groups}
Irreducible sextics~$C$ with a singular point of
multiplicity~$4$ and nonabelian fundamental group form seven
rigid isotopy classes, one class for each of the following
formulas\rom:
\dashes
\dash
$\{(2,0),(2,0),(0,5),(0,5)\}$, with
$\pi_1(\Cp2\sminus C)=G_2\cong\GD_{10}\times(\ZZ/3\ZZ)$\rom;
\dash
$\{(4,0),(0,5),(0,5)\}$, with $\pi_1(\Cp2\sminus C)=G_4$\rom;
\dash
$\{(2,2),(2,2),(0,3),(0,3)\}$\rom;
\dash
$\{(2,5),(2,2),(0,3)\}$\rom;
\dash
$\{(2,5),(2,5)\}$\rom;
\dash
$\{(4,4),(0,3),(0,3)\}$\rom;
\dash
$\{(4,7),(0,3)\}$.
\enddashes
In the first two cases, the curve is not of torus type and one has
$\Delta_C(t)=1$\rom;
in the last five cases, the curve is of torus type and one has
$\pi_1(\Cp2\sminus C)=\GB_3/\Delta^2$.
\endtheorem

\proof
The fundamental group of a curve~$C$ of degree~$m$ with a singular
point of multiplicity~$m-2$ is described in~\cite{groups}. If
$m=6$, it is easy to enumerate all formulas satisfying
conditions~\itemref{s.m4}{m4.1}--\ditto{m4.3} above and select those that
give rise to non-abelian fundamental groups. Then,
Theorem~\ref{m-2.classes} would apply to give one rigid isotopy
class for each formula found.

In fact, for $\pi_1(\Cp2\sminus C)$ not to be abelian, one must
have $p_i\ne1$ for all~$i$; since also $\sum p_i=4$,
the nonzero entries~$p_i$ are either $2$ and~$2$ or~$4$. Another
necessary condition is that $q=\gcd\{q_i-p_i\}_{q_i>p_i}$ must be
larger than~$1$ and, in view of~\ditto{m4.3} above, $q$ must also be
odd. These restrictions leave the seven formulas listed in the
statement.

For the last five formulas one has
$q_i\ge p_i$ and $3\mathbin|(q_i-p_i)$ for each~$i$. Hence, $a$~is
a square, $a\mathbin|D$, and $D/a$ is a cube. Let $a=p^2$ and
$D=-4p^2s^3$. Then $p\mathbin|b$, $b=2pq$, and $c=q^2+s^3$. Thus,
the equation of~$C$ has the form
$$
(x_0p+q)^2+s^3=0,
$$
\ie, $C$ is of torus type.
\endproof


\let\ref\bibitem
\let\endref\relax
\let\publ\relax
\let\publaddr\relax
\def\by#1{{\bibname #1},}
\def\paper#1{`#1',}
\def\book#1{`#1',}
\def\jour#1{{\em #1\/}}
\def\yr#1{(#1)}
\def\issue#1{no.~#1,}
\def\vol#1{{#1}}
\def\pages#1{#1.}
\def\finalinfo#1{{#1}}

\affiliationone{
Alex Degtyarev\\
Department of Mathematics,
Bilkent University,
06800 Ankara\\
Turkey
\email{degt\at fen.bilkent.edu.tr}}

\end{document}